\input ./style/arxiv-general.cfg
\documentclass[aop,MSNbibl,dvips]{arximspdf}
\makeatletter
   \@ifpackageloaded{graphicx}{}{\usepackage{graphicx}}
\makeatother
\usepackage{mathbh}

\doi{10.1214/14-AOP952} 
\volume{43}
\issue{5}
\pubyear{2015}
\firstpage{2810}
\lastpage{2840}
\docsubty{FLA}

\makeatletter
\newcommand{\rrvert}{\vert}
\newcommand{\rrVert}{\Vert}
\newcommand{\llvert}{\vert}
\newcommand{\llVert}{\Vert}
\renewcommand{\mid}{|}
\newtheorem{introthm}{Theorem}
\newtheorem{theorem}{Theorem}[section]
\newtheorem{proposition}[theorem]{Proposition}
\newtheorem{lemma}[theorem]{Lemma}
\newtheorem{conjecture}[theorem]{Conjecture}
\newproclaim{remark}[theorem]{Remark}
\newproclaim{notation}[theorem]{Notation}
\newproclaim{definition}[theorem]{Definition}
\newproclaim{example}[theorem]{Example}
\newcommand{\Aut}{\operatorname{Aut}}
\newcommand{\St}{\operatorname{Stab}}
\newcommand{\var}{\operatorname{var}}
\newcommand{\cov}{\operatorname{cov}}
\newcommand{\corr}{\operatorname{corr}}
\newcommand{\vol}{\operatorname{vol}}
\newcommand{\spn}{\operatorname{span}}
\newcommand{\cl}{\operatorname{cl}}
\newcommand{\dist}{\operatorname{dist}}
\def\ind{\mathbh{1}} 
\makeatother

\begin{document}
\begin{frontmatter}

\title{Independence ratio and random eigenvectors in~transitive graphs}
\runtitle{Independence ratio and random eigenvectors}

\begin{aug}
\author{\fnms{Viktor}~\snm{Harangi}\corref{}\thanksref{T1}\ead[label=e1]{harangi@math.toronto.edu}}
\and
\author{\fnms{B\'alint}~\snm{Vir\'ag}\thanksref{T2}\ead[label=e2]{balint@math.toronto.edu}}
\runauthor{V. Harangi and B. Vir\'ag}
\affiliation{University of Toronto}
\address{Department of Mathematics\\
University of Toronto\\
40 St. George Street\\
Toronto, Ontario M5S 2E4\\
Canada\\
\printead{e1}\\
\phantom{E-mail: }\printead*{e2}} 
\end{aug}
\thankstext{T1}{Supported by MTA R\'enyi ``Lend\"ulet'' Groups and
Graphs Research Group.}
\thankstext{T2}{Supported by the Canada Research Chair and NSERC DAS programs.}

\received{\smonth{8} \syear{2013}}
\revised{\smonth{6} \syear{2014}}

%
\begin{abstract}
A theorem of Hoffman gives an upper bound
on the independence ratio of regular graphs
in terms of the minimum $\lambda_{\min}$ of the spectrum of the
adjacency matrix.
To complement this result we use random eigenvectors
to gain lower bounds in the vertex-transitive case.
For example, we prove that the independence ratio
of a $3$-regular transitive graph is at least
\[
q = \frac{1}{2} - \frac{3}{4 \pi} \arccos\biggl( \frac{1-\lambda
_{\min}}{4}
\biggr).
\]
The same bound holds for infinite transitive graphs:
we construct factor of i.i.d. independent sets for which
the probability that any given vertex is in the set is at least $q-o(1)$.

We also show that the set of the distributions of factor of i.i.d. processes
is not closed w.r.t. the weak topology
provided that the spectrum of the graph is uncountable.
\end{abstract}

%
\begin{keyword}[class=AMS]
\kwd{05C69}
\kwd{05C50}
\kwd{60G15}
\end{keyword}
\begin{keyword}
\kwd{Independent set}
\kwd{independence ratio}
\kwd{minimum eigenvalue}
\kwd{adjacency matrix}
\kwd{regular graph}
\kwd{transitive graph}
\kwd{factor of i.i.d.}
\kwd{invariant Gaussian process}
\end{keyword}
\end{frontmatter}

\setcounter{footnote}{2}


\section{Introduction}\label{sec1}

\subsection{The independence ratio and the minimum eigenvalue}

An \emph{independent set} is a set of vertices in a graph, no two of
which are adjacent.
The \emph{independence ratio} of a graph $G$ is the size of
its largest independent set divided by the total number of vertices.
If $G$ is regular, then the independence ratio is at most $1/2$,
and it is equal to $1/2$ if and only if $G$ is bipartite.

The adjacency matrix of a $d$-regular graph has
real eigenvalues between $-d$~and~$d$.
The least eigenvalue $\lambda_{\min}$ is at least $-d$,
and it is equal to $-d$ if and only if the graph is bipartite.

So the distance of the independence ratio from $1/2$
and the distance of $\lambda_{\min}$ from $-d$ both measure
how far a $d$-regular graph is from being bipartite.
The following natural question arises: what kind of connection
is there between these two graph parameters?

A theorem of Hoffman \cite{hoffman1} gives a partial answer to this question.
It says that the independence ratio of a $d$-regular graph is at most
%
%
\begin{equation}
\label{eqhoffman} \frac{-\lambda_{\min}}{d-\lambda_{\min}} = \frac
{1}{2} - \frac{(1/2)(\lambda_{\min}+ d)}{2d - (
\lambda_{\min}+ d )};
\end{equation}
for a simple proof, see \cite{ellis1}, Theorem~11;
also see \cite{lyonsnazarov}, Section~4, for certain improvements.

Hoffman's bound implies that $\lambda_{\min}\to-d$ as the
independence ratio tends to~$1/2$.
The converse statement is not true in general:
it is easy to construct \mbox{$d$-}regular graphs with $\lambda_{\min}$
arbitrarily close to $-d$
and the independence ratio separated from $1/2$.
However, for transitive graphs the converse is also true.
A graph $G$ is said to be \emph{vertex-transitive} (or \emph
{transitive} in short)
if its automorphism group $\Aut(G)$ acts transitively on the vertex
set $V(G)$.
%

\begin{introthm} \label{thmcrude}
Let $G$ be a finite, $d$-regular, vertex-transitive graph with
least eigenvalue $\lambda_{\min}$.
Then the independence ratio of $G$ is at least
\[
\tfrac{1}{2} - \tfrac{1}{3} \sqrt{ d(\lambda_{\min}+ d) }.
\]
In particular, if $\lambda_{\min}\to-d$, then the independence ratio
converges to $1/2$.
\end{introthm}

The idea behind the proof is to consider random eigenvectors with
eigenvalue $\lambda_{\min}$.
Let $\lambda$ be an arbitrary eigenvalue of the adjacency matrix
of some transitive graph $G$, and let $E_\lambda$ denote the
eigenspace corresponding to $\lambda$,
that is, the space of eigenvectors with eigenvalue $\lambda$.
(Note that $E_\lambda$ is typically more than one dimensional, since
$G$ is transitive.)
Furthermore, let $S_\lambda$ be the unit sphere in $E_\lambda$.
Now we pick a uniform random vector from $S_\lambda$.
Note that $S_\lambda$ is $\Aut(G)$-invariant, therefore
the distribution of this random vector is $\Aut(G)$-invariant, too.
Let us choose the vertices $v$ with the property that
the value of the eigenvector at $v$ is larger than at each neighbor of $v$.
(If $\lambda$ is negative, then
we expect many of the vertices with positive value to have this property.)
Clearly, these vertices form an independent set.
Since our random vector is invariant, the probability $q$
that a given vertex is chosen is the same for all vertices.
Therefore the expected size of this random independent set is $q
\llvert V(G)\rrvert $,
and consequently, the independence ratio of $G$ is at least $q$.
An estimate of $q$ yields Theorem~\ref{thmcrude} above.
In many cases we obtain much sharper bounds.

When the graph has a lot of symmetry (e.g.,
when any pair of neighbors of a fixed vertex can be mapped
to any other pair by a suitable graph automorphism),
then the probability $q$ defined above is actually determined by
$\lambda$.
In this case it equals $q_d(\lambda)$, the relative volume of
the $(d-1)$-dimensional regular spherical simplex defined by normal vectors
with pairwise scalar product $\frac{d-2-\lambda}{2(d-1)}$; see
Definition~\ref{defqd}.
There is a simple formula for $q_3(\lambda)$; see Theorem~\ref
{thm3regmain}.

We conjecture that $q \geq q_d(\lambda)$ for arbitrary transitive graphs
(provided that $\lambda$ is sufficiently small).
In other words, the worst-case scenario is when the graph has a lot of symmetry.
Of course, this would yield a lower bound $q_d(\lambda_{\min})$ for
the independence ratio.
We managed to prove this conjecture for $3$-regular transitive graphs and
$4$-regular arc-transitive graphs. We also showed that
a natural geometric conjecture would imply the $d$-regular,
arc-transitive case.
[A graph is said to be \emph{arc-transitive} or \emph{symmetric}
if for any two pairs of adjacent vertices $(u_1,v_1)$ and $(u_2,v_2)$,
there is an automorphism of the graph mapping $u_1$ to $u_2$ and $v_1$
to $v_2$.]
The following theorems were obtained.

%
%
\begin{introthm} \label{thmarctr}
Suppose that $G$ is a finite, $d$-regular, arc-transitive graph
with least eigenvalue $\lambda_{\min}$.
Then the independence ratio of $G$ is at least
\[
\tfrac{1}{2} - \tfrac{1}{3} \sqrt{\lambda_{\min}+d}.
\]
In fact, a natural geometric conjecture (see Conjecture~\ref{conjgeom})
would imply that the independence ratio is at least $q_d(\lambda_{\min})$.
This has been proven in the case $d=4$: the independence ratio of
a finite, $4$-regular, arc-transitive graph is at least
%
%
\begin{equation}
\label{eqq4} q_4(\lambda_{\min}) \geq\tfrac{1}{2} -
\tfrac{1}{4} \sqrt{\lambda_{\min}+4}.
\end{equation}
\end{introthm}

%
\begin{introthm} \label{thm3regmain}
Suppose that $G$ is a finite, $3$-regular, vertex-transitive graph
with minimum eigenvalue $\lambda_{\min}$. Then the independence ratio
of $G$ is at least
\[
q_3(\lambda_{\min}) = \frac{1}{8} + \frac{3}{4 \pi}
\arcsin\biggl( \frac{1-\lambda_{\min}}{4} \biggr) = \frac{1}{2} -
\frac{3}{4 \pi} \arccos\biggl( \frac{1-\lambda
_{\min}}{4} \biggr).
\]
In fact, the following stronger statement holds:
$G$ contains two disjoint independent sets $I_1, I_2$
with total size $\llvert I_1 \cup I_2\rrvert \geq2 q_3(\lambda_{\min
}) \llvert V(G)\rrvert $.
This means that the induced subgraph $G[I_1 \cup I_2]$ is bipartite
and has at least $2 q_3(\lambda_{\min}) \llvert V(G)\rrvert $ vertices.
\end{introthm}

%
%
\begin{figure}[t]

\includegraphics{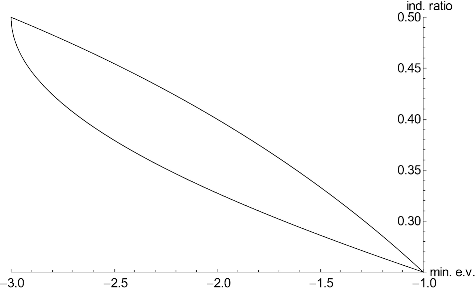}

\caption{Hoffman's upper bound (\protect\ref{eqhoffman})
and the lower bound of Theorem \protect\ref{thm3regmain} for
$\lambda_{\min}\in[-3,-1]$.}
\label{figcompare}
\end{figure}

See Figure~\ref{figcompare} to compare the lower bound given
in Theorem~\ref{thm3regmain} to Hoffman's upper bound (\ref{eqhoffman}).
Note that $-3 \leq\lambda_{\min}\leq-2$ for any $3$-regular
transitive graph
with the only exception of the complete graph $K_4$ for which $\lambda
_{\min}= -1$; see Proposition~\ref{proplamin3reg} in the \hyperref[secapp]{Appendix}.

\subsection{Random wave functions on infinite transitive graphs}

In order to generalize the above theorems,
we define random wave functions on infinite transitive graphs $G$.
A \emph{wave function} with eigenvalue $\lambda$ on $G$ is a function
$f \colon V(G) \to\mathbb{R}$ such that
\[
\sum_{u \in N(v)} f(u) = \lambda f(v)\qquad\mbox{for each
vertex } v \in V(G),
\]
where $N(v)$ denotes the set of neighbors of $v$ in $G$.
So a wave function is basically an eigenvector of the adjacency
operator of $G$,
except that it does not need to be in $\ell_2(V(G))$.

These random wave functions will also let us answer
an open question concerning factor of i.i.d. processes.
Suppose that we have independent standard normal random variables $Z_u$
assigned to each vertex $u$ of an infinite transitive graph $G$.
By a \emph{factor of i.i.d. process} on $G$ we mean random variables
$X_v$, $v \in V(G)$
that are all obtained as measurable functions of the random variables
$Z_u$, $u \in V(G)$
and that are $\Aut(G)$-equivariant [i.e., they commute with the
natural action of $\Aut(G)$].
It is easy to see that for any factor of i.i.d. process $X_v$, $v \in V(G)$
with $0< \var(X_v) < \infty$, the correlation of $X_v$ and $X_{v'}$
converges to $0$
as the distance of $v$ and $v'$ goes to infinity; see Proposition~\ref
{propcorrdecay} in the \hyperref[secapp]{Appendix}.
So a random process that is $0$ everywhere with probability $1/2$
and $1$ everywhere with probability $1/2$ cannot be a factor of i.i.d.
However, it can be seen easily that this process can be approximated
by factor of i.i.d. processes provided that $G$ is amenable.
So the space of factor of i.i.d. processes is not closed; that is,
the distributions of these processes do not form a closed set w.r.t.
the weak topology.
It has been an open question
whether the same is true on nonamenable graphs,
for example, on the \mbox{$d$-}regular tree; see \cite{birsreport}, Section~4, Question~4.
We will show that the space of factor of i.i.d. processes
is not closed provided that the spectrum of $G$ is uncountable.

We say that a factor of i.i.d. process $X_v$, $v \in V(G)$
is a \emph{linear factor of i.i.d.} if each $X_v$ is
obtained as a (possibly infinite) linear combination of $Z_u$, $u \in V(G)$.
Note that linear factors have the following properties.
%

\begin{definition} \label{defgaussianprocess}
We call a collection of random variables $X_v$, $v \in V(G)$
a~\emph{Gaussian process} on $G$ if
they are jointly Gaussian, and each $X_v$ is centered (i.e., has mean $0$).
(Random variables are jointly Gaussian
if any finite linear combination of them is Gaussian.)
We say that a Gaussian process $X_v$ is $\Aut(G)$-invariant (or simply
invariant) if
for any $\Phi\in\Aut(G)$ the joint distribution of
the Gaussian process $X_{\Phi(v)}$ is the same as that of the original process.
\end{definition}

We will prove that the adjacency operator $A_G$ has
approximate eigenvectors (satisfying a certain invariance property)
for any $\lambda$ in the spectrum $\lambda\in\sigma(A_G)$.
Then we will use these approximate eigenvectors
as coefficients to define linear factor of i.i.d. processes
converging in distribution to an invariant Gaussian process $X_v$
that satisfies the eigenvector equation at each vertex.
%

\begin{introthm} \label{thmgaussianev}
Let $G$ be an infinite vertex-transitive graph with adjacency operator $A_G$.
Then for each point $\lambda$ of the spectrum $\sigma(A_G)$ there exists
a nontrivial invariant Gaussian process $X_v$, $v \in V(G)$ such that
%
%
\begin{equation}
\label{eqeigen} \sum_{u \in N(v)} X_u = \lambda
X_v\qquad\mbox{for each vertex } v \in V(G),
\end{equation}
where $N(v)$ denotes the set of neighbors of $v$ in $G$.
Furthermore, the process $X_v$ can be approximated (in distribution)
by linear factor of i.i.d. processes. Clearly, we can assume that
these approximating linear factors have only finitely many nonzero coefficients.
\end{introthm}

An invariant Gaussian process satisfying (\ref{eqeigen}) will be called
a \emph{Gaussian wave function} with eigenvalue $\lambda$.
If the spectrum of $G$ is not countable, then
we can conclude that some of these Gaussian wave functions
cannot be obtained as factor of i.i.d. processes.
%

\begin{introthm} \label{thmnotclosed}
Let $G$ be an infinite transitive graph such that
the spectrum of the adjacency oparator $A_G$ is not countable.
Then there exist (linear) factor of i.i.d. processes on $G$ with the property
that the weak limit of their distributions cannot be
obtained as the distribution of a factor of i.i.d. process.
\end{introthm}

We can say more for Cayley graphs.
%

\begin{introthm} \label{thmcayley}
Suppose that $G$ is the Cayley graph of a finitely generated infinite group.
Then a Gaussian wave function with eigenvalue $\lambda_{\max
}\stackrel{\mathrm{def}}{=}
\sup\sigma(A_G)$
can never be obtained as the distribution of a factor of i.i.d. process.
\end{introthm}

In view of Theorems~\ref{thmgaussianev} and~\ref{thmcayley}
there exists a Gaussian wave function with eigenvalue $\lambda_{\max
}$ that
can be approximated by factor of i.i.d. processes but cannot be
obtained as one.
An independent and different proof of this result
was given by Russell Lyons in the special case
when $G$ is a regular tree \cite{lyons}, Corollary~3.3.

\subsection{Factor of i.i.d. independent sets}

Let $X_v$, $v \in V(G)$ be a random process
on our infinite transitive graph $G$.
As in the finite setting,
$I_{+} \stackrel{\mathrm{def}}{=} \{ v\dvtx  X_v > X_u$,
$\forall u \in N(v)
\}$
is a random independent set.
If our process is invariant, then the probability that $v \in I_+$
is the same for each vertex $v$, and thus this probability
can be used to measure the size of $I_+$.
If our process is a factor of some i.i.d. process $Z_v$,
then the resulting independent set is also a factor of $Z_v$.

In the infinite setting let $\lambda_{\min}$ denote the minimum of
the spectrum $\sigma(A_G)$,
and let $X_v$ be a linear factor of $Z_v$ approximating
the Gaussian eigenvector with eigenvalue $\lambda_{\min}$; see
Theorem~\ref{thmgaussianev}.
As the process $X_v$ converges in distribution to the Gaussian eigenvector,
the probability $P( v \in I_+ )$ approaches
the corresponding probability for the Gaussian eigenvector process,
which, as we will see, can be computed the exact same way
as in the finite case.
%

\begin{introthm} \label{thminf}
Theorems~\ref{thmcrude},~\ref{thmarctr} and~\ref{thm3regmain}
give lower bounds $q$ (in terms of $\lambda_{\min}$) for the
independence ratio
of finite transitive graphs with least eigenvalue $\lambda_{\min}$.
These bounds remain true in the following framework.
Let $\lambda_{\min}$ denote the minimum of the spectrum of an
infinite transitive graph $G$.
Then for any $\varepsilon> 0$ there exists a factor of i.i.d.
independent set on $G$
such that the probability that any given vertex is in the set is at
least $q-\varepsilon$.
\end{introthm}

A special case of this infinite setting was investigated in \cite{csghv}.
When $G$ is the $d$-regular tree $T_d$,
then any factor of i.i.d. independent set on $G$
automatically gives a lower bound for the independence ratio of
$d$-regular finite graphs with sufficiently large girth.
In particular, for the $3$-regular tree, $T_3$ one has $\lambda_{\min
}= -2\sqrt{2}$.
Therefore the infinite version of Theorem~\ref{thm3regmain} tells us that
there exists factor of i.i.d. independent set in $T_3$ with density
\[
\frac{1}{2} - \frac{3}{4\pi} \arccos\biggl( \frac{1+2\sqrt{2}}{4}
\biggr) \approx0.4298.
\]
In \cite{csghv} the somewhat better bound $0.4361$ was obtained.
In fact, \cite{csghv} was the starting point for the work in the
present paper.
For previous results on the independence ratio of large-girth graphs,
see \cite{Boindset,mckay,shearer,shearer2,lauerwormald,cubic}.


\section{Finite vertex-transitive graphs} \label{sec2}

Throughout this section $G$ will denote a vertex-transitive, finite graph
with degree $d$ for some positive integer $d \geq3$.
The least eigenvalue of its adjacency matrix $A_G$ will be denoted by
$\lambda_{\min}$.
For now let $\lambda$ be an arbitrary eigenvalue of $A_G$.
Eventually, we will choose $\lambda$ as the minimum eigenvalue.
First we define what we mean by a random eigenvector.
%

\begin{definition} \label{defrandomev}
Let $E_\lambda$ be the eigenspace corresponding to $\lambda$, that is,
\[
E_\lambda\stackrel{\mathrm{def}} {=} \bigl\{ x \in\ell_2
\bigl( V(G) \bigr)\dvtx  A_G x = \lambda x \bigr\}.
\]
We fix some orthonormal basis $e_1, \ldots, e_l$ in $E_\lambda$,
and take independent standard normal random variables $\gamma_1,\ldots, \gamma_l$.
We call $\sum_{i=1}^l \gamma_i e_i$
the \emph{random eigenvector with eigenvalue $\lambda$}.
\end{definition}
%

\begin{remark}
The (distribution of the) random eigenvector is clearly independent of
the choice of the basis $e_1, \ldots, e_l$, so it is well defined.
It also follows that the distribution of the random eigenvector is
$\Aut(G)$-invariant.
(Note that in the \hyperref[sec1]{Introduction} we defined the random eigenvector differently:
a uniform random vector on the unit sphere of $E_\lambda$,
which is just the normalized version of the random eigenvector of
Definition~\ref{defrandomev}.)
\end{remark}

We will think of this random eigenvector as a collection of
real-valued random variables $X_v$, $v \in V(G)$ with the property that
they are jointly Gaussian and $\Aut(G)$-invariant, each $X_v$ is
centered, and
\[
\sum_{u \in N(v)} X_u = \lambda
X_v\qquad\mbox{for each vertex } v, %
\]
where $N(v)$ denotes the set of neighbors of $v$ in $G$.
Since $G$ is transitive, each $X_v$ has the same variance.
After multiplying these random variables with a suitable positive constant
we might assume that $\var(X_v) = 1$ for each vertex $v$.
Next we define random independent sets by means of these random eigenvectors.
%

\begin{definition} \label{defindset}
Let $X_v$, $v \in V(G)$ denote the random eigenvector corresponding to
the eigenvalue $\lambda$
as explained above. The random sets $I_{+}$ and $I_{-}$ are defined as follows:
\begin{eqnarray*}
I_{+} &=& I^\lambda_{+} \stackrel{\mathrm{def}} {=}
\bigl\{ v\in V(G)\dvtx  X_v > X_u\mbox{ for each } u \in
N(v) \bigr\}\quad\mbox{and}
\\
I_{-} &=& I^\lambda_{-} \stackrel{\mathrm{def}} {=}
\bigl\{ v\in V(G)\dvtx  X_v < X_u\mbox{ for each } u \in
N(v) \bigr\}.
\end{eqnarray*}
Clearly, $I_{+}$ and $I_{-}$ are disjoint (random) independent sets in $G$.
\end{definition}

The $\Aut(G)$-invariance implies that the probability of
the event $v \in I_{+}$ is the same for all vertices $v$.
So from now on, we will focus on a fixed vertex and its neighbors.
First we introduce the following notation.
%

\begin{notation} \label{notsetup}
Let $v$ be an arbitrary vertex of our vertex-transitive graph $G$.
We will call $v$ the root. The neighbors of $v$ are denoted by $w_1,
\ldots, w_d$.
For $X_v$ and $X_{w_i}$ we will simply write $X$ and $Y_i$, respectively.
We will assume that $\var(X) = 1$, which implies that $\var(Y_i) = 1$
for each $i$.
Since $X_v$, $v \in V(G)$ is the random eigenvector with eigenvalue
$\lambda$, we have
%
%
\begin{equation}
\label{eqev} \sum_{i=1}^d
Y_i = \lambda X.
\end{equation}
The covariance $\cov(Y_i,Y_j)$ will be denoted by $c_{i,j}$.
Then it readily follows from~(\ref{eqev}) that
%
%
\begin{equation}
\label{eqcovsum}\lambda^2 = \cov(\lambda X, \lambda X) = \sum
_{i,j} c_{i,j} = d + 2 \sum
_{i<j} c_{i,j},\qquad\mbox{thus } \sum
_{i<j} c_{i,j} = \frac{\lambda^2 - d}{2}.\hspace*{-20pt}
\end{equation}
We also introduce the random variables $ U_i \stackrel{\mathrm
{def}}{=}X - Y_i$.

Since $Y_1, \ldots, Y_d$ are centered and jointly Gaussian, they can be
written as the linear combinations of independent standard normal variables:
there exist independent standard Gaussians $Z_1, \ldots, Z_d$
and (deterministic) vectors $y_1, \ldots, y_d \in\mathbb{R}^d$ such that
$Y_i$ is the inner product of $y_i$ and $Z=(Z_1,\ldots,Z_d)$.
Setting $x = (y_1 + \cdots+ y_d)/\lambda$ and $u_i = x - y_i$, we have
\[
Y_i = y_i \cdot Z;\qquad
X = x \cdot Z;\qquad
U_i = u_i \cdot Z.
\]
It is easy to see that for any deterministic vectors $a,b \in\mathbb{R}^d$
the covariance\break $\cov(a \cdot Z, b \cdot Z)$ is equal to the inner
product $a \cdot b$.
In particular,
%
%
\begin{eqnarray}\label{eqcovinnerpr}
x \cdot x &=& \var(X) = 1;\qquad
y_i \cdot y_j = \cov(Y_i, Y_j) = c_{i,j};
\nonumber\\[-8pt]\\[-8pt]\nonumber
u_i \cdot u_j &=& \cov(U_i, U_j).
\end{eqnarray}
Finally, we introduce the following notation for the pairwise angles of
the vectors~$u_i$:
%
%
\begin{equation}
\label{eqphi} \varphi_{i,j} \stackrel{\mathrm{def}} {=}
\angle(u_i,u_j) = \arccos\biggl( \frac{u_i \cdot u_j}{\llVert u_i
\rrVert \llVert u_j \rrVert }
\biggr).
\end{equation}
\end{notation}

Our goal is to give estimates for the probability that
a certain vertex lies in our random independent set $I_{+}$.
As we will see, this probability can be expressed as the volume of a
certain spherical simplex.
%

\begin{definition}
Let $S^{d-1}$ denote the unit sphere in $\mathbb{R}^d$.
A half-space is said to be \emph{homogeneous} if the defining hyperplane
(i.e., the boundary of the half-space) passes through the origin.
A vector $n$ orthogonal to the defining hyperplane and ``pointing outward''
is called an \emph{outer normal vector}.
Then the given (open) half-space consists of those $x \in\mathbb{R}^d$
for which the inner product $n \cdot x$ is negative.\vspace*{1pt}

A $(d-1)$-dimensional \emph{spherical simplex} is
the intersection of $S^{d-1}$ and $d$ homogeneous half-spaces in
$\mathbb{R}^d$.
Up to congruence, a spherical simplex is determined by the ${d\choose2}$
pairwise angles enclosed by the outer normal vectors of the $d$ half-spaces.
If these ${d\choose2}$ angles are all equal, then we say that
the spherical simplex is \emph{regular}.
\end{definition}
%

\begin{proposition} \label{propsimplexprob}
The probability that any fixed vertex is in
the random independent set $I_{+}$
is equal to the \emph{relative volume} of the $(d-1)$-dimensional
spherical simplex
corresponding to the outer normal vectors $-u_i$.
\end{proposition}

\begin{pf}
The probability in question is
%
%
\begin{equation}
\label{eqorthant} P(v \in I_{+}) = P( X > Y_i, 1 \leq i
\leq d ) = P( U_i > 0, 1 \leq i \leq d ).
\end{equation}
The event $U_i > 0$ is that the random point $Z$ lies
in the homogeneous open half-space with outer normal vector $-u_i$.
So the probability is equal to the measure of
the intersection of the homogeneous half-spaces with outer normal
vectors $-u_i$
with respect to the standard multivariate Gaussian measure on $\mathbb{R}^d$.
This is simply the volume of the corresponding $(d-1)$-dimensional
spherical simplex
divided by the volume $\vol(S^{d-1})$ of the unit sphere $S^{d-1}$.
Note that this relative volume is determined by the pairwise angles
$\varphi_{i,j}$
[see (\ref{eqphi}) in Notation~\ref{notsetup}],
which, in turn, can be expressed using the inner products $y_i \cdot
y_j = c_{i,j}$.
\end{pf}

The probability $P(v \in I_{+})$ seems to be the smallest
when $G$ has a lot of symmetry.
To make this more precise, we first define what we mean by a ``lot of
symmetry.''
%

\begin{definition}
We say that $G$ is \emph{cherry-transitive} if
any cherry (path of length $2$) in $G$ can be mapped to any other cherry
using a suitable graph automorphism of $G$.
\end{definition}
%

\begin{proposition} \label{propchtr}
If $G$ is cherry-transitive, then
\[
c_{i,j} = \frac{\lambda^2 - d}{d(d-1)}\qquad \mbox{for all } i \neq j,
\]
and, consequently, the pairwise angles $\varphi_{i,j}$ are all equal to
%
%
\begin{equation}
\label{eqangle} \arccos\biggl( \frac{d-2-\lambda}{2(d-1)} \biggr).
\end{equation}
\end{proposition}

\begin{pf} 
If $G$ is cherry-transitive, then for any $i_1 \neq j_1$ and $i_2 \neq
j_2$ there exists
an automorphism $\Phi\in\Aut(G)$ such that $\Phi$ fixes the root
$v$ and
takes the unordered pair $w_{i_1}, w_{j_1}$ to $w_{i_2}, w_{j_2}$, that is,
\[
\Phi v = v,\qquad
\Phi w_{i_1} = \Phi w_{i_2},\qquad
\Phi w_{j_1} = \Phi w_{j_2}
\]
or
\[
\Phi v = v,\qquad
\Phi w_{i_1} = \Phi w_{j_2},\qquad
\Phi w_{j_1} = \Phi w_{i_2}.
\]
Together with the $\Aut(G)$-invariance of the random eigenvector
this implies that $c_{i_1,j_1} = c_{i_2,j_2}$.
Since this holds for any two pairs of indices,
it follows that all $c_{i,j}$, $i \neq j$ are the same.
Using (\ref{eqcovsum}) we conclude that for $i \neq j$
\[
c_{i,j} = \frac{\lambda^2 - d}{d(d-1)}.
\]
Then easy calculation shows
(using notations introduced earlier) that
\[
\llVert u_i \rrVert^2 = \llVert u_j
\rrVert^2 = \frac{2(d-\lambda)}{d}\quad\mbox{and}\quad u_i \cdot
u_j = \frac{(d-\lambda)(d-2-\lambda)}{d(d-1)}.
\]
Plugging this into (\ref{eqphi}) gives
\[
\varphi_{i,j} = \arccos\biggl( \frac{d-2-\lambda}{2(d-1)} \biggr).
\]\upqed
\end{pf}

We are now in a position to define the functions $q_d(\lambda)$.
%

\begin{definition} \label{defqd}
For $-d \leq\lambda\leq d$, let $q_d(\lambda)$ denote the volume of
the $(d-1)$-dimensional regular spherical simplex corresponding to
the angle (\ref{eqangle}) divided by $\vol(S^{d-1})$.
Then $P(v \in I_{+}) = q_d(\lambda)$ for any cherry-transitive $G$.
In particular, the independence ratio of
any cherry-transitive graph $G$ is at least $q_d(\lambda_{\min})$.
\end{definition}

So $P(v \in I_{+}) = q_d(\lambda)$ provided that $G$ has enough symmetry.
The following conjecture says that in the general (i.e.,
vertex-transitive) case
the probability should be larger than that.
%

\begin{conjecture} \label{conjqd}
For any transitive graph $G$ it holds that
\[
P(v \in I_{+}) \geq q_d(\lambda) %
\]
for any $\lambda$, or at least for sufficiently small $\lambda$:
$\lambda\leq\lambda_0$ for some $\lambda_0$.

This would, of course, imply that the independence ratio of $G$ is at
least $q_d(\lambda_{\min})$
provided that $\lambda_{\min}\leq\lambda_0$.
\end{conjecture}

We will prove this conjecture for $d=3$ and $\lambda_0=-2$ in
Section~\ref{sec3reg}.
The conjecture might be true for arbitrary $\lambda$,
but proving for $\lambda\leq\lambda_0 = -2$ will be sufficient for
our purposes,
because $\lambda_{\min}\leq-2$ for any $3$-regular transitive graph
except $K_4$.
%

\begin{remark}
In view of formula (\ref{eqorthant})
the above conjecture would follow from the following statement.
Let $(U_1, \ldots, U_d)$ be a multivariate Gaussian with each $U_i$ centered
and with pairwise covariances $m_{i,j}$.
We have the following constraints for the covariances:
\[
\sum_{1 \leq i,j \leq d} m_{i,j} = (d-
\lambda)^2\qquad\mbox{and for each } i\dvtx   m_{i,i} =
\frac{2}{d-\lambda} \sum_{j=1}^{d}
m_{i,j}.
\]
[Or\vspace*{1pt} one might replace these by the weaker constraints
$\sum_{1 \leq i,j \leq d} m_{i,j} = (d-\lambda)^2$ and $\sum_{1
\leq i \leq d} m_{i,i} = 2(d-\lambda)$.]\vspace*{1pt}
Then the orthant probability $P( U_i > 0, 1 \leq i \leq d )$
is minimized when all $m_{i,i}$, $1\leq i \leq d$ are equal
and also all $m_{i,j}$, $i \neq j$ are equal.
\end{remark}

A few properties of the functions $q_d(\lambda)$ are collected in the
next proposition.
%

\begin{proposition} \label{propqdprop}
For any $d \geq3$, $q_d$ is a monotone decreasing continuous function
on $[-d,-1]$ with
\[
q_d(-d) = \frac{1}{2}\quad\mbox{and}\quad q_d(-1) =
\frac{1}{d+1}.
\]
As for the behavior of $q_d$ around $-d$ we have
\[
q_d(\lambda) \geq\frac{1}{2} - \frac{ \pi\vol(S^{d-2}) }{ 4 \vol
(S^{d-1}) } \sqrt{
\frac{\lambda+d}{d}} \geq\frac{1}{2} - \frac{1}{3} \sqrt{\lambda+d}.
\]
\end{proposition}

\begin{pf}
The volume of a regular spherical simplex is clearly a continuous and
monotone decreasing
function of the corresponding angle. Since (\ref{eqangle}) is
a continuous and monotone increasing function of $\lambda\in[-d,-1]$,
monotonicity and continuity of $q_d$ follow.

For $\lambda= -d$ the angles $\varphi_{i,j}$ are $0$, so
the corresponding (degenerate) spherical simplex is a hemisphere,
thus $q_d(-d) = 1/2$ as claimed.

For $\lambda= -1$ the angles $\varphi_{i,j}$ are $\pi/3$.
It is not hard to see that the vertices of our spherical simplex in
that case
will be the $d$ vertices of a face of a regular (Euclidean) simplex in
$\mathbb{R}^d$.
Then each of the $d+1$ spherical simplices belonging
to the $d+1$ faces has volume $\vol(S^{d-1})/(d+1)$.
[We could also argue that for $G = K_{d+1}$ and $\lambda= -1$
we have $P( v \in I_{+} ) = 1/(d+1)$,
and since $K_{d+1}$ is cherry-transitive,
$P( v \in I_{+} ) = q_d(-1)$.]

See Section~\ref{secnearnegd} for a proof of the claimed behavior
around $-d$.
\end{pf}

\subsection{The 3-regular, vertex-transitive case} \label{sec3reg}

Now we turn to the proof of Theorem~\ref{thm3regmain} that gives
a lower bound for the independence ratio of $3$-regular transitive graphs.
We will basically show that Conjecture~\ref{conjqd}
is true when $d=3$ and $\lambda_0 = -2$.

As we have seen, $P(v \in I_{+})$ equals
the relative volume of a certain spherical simplex.
For $d=3$ the surface of the unit sphere $S^{d-1} = S^2$ is $4\pi$ and
the area of a spherical triangle with interior angles $\alpha, \beta, \gamma$
is equal to $\alpha+ \beta+ \gamma- \pi$.
The spherical triangle in question is determined by the homogeneous half-spaces
with outer normal vectors $-u_1, -u_2, -u_3$ (recall Proposition~\ref
{propsimplexprob}).
We denoted the angle enclosed by the outer normal vectors $-u_i$ and
$-u_j$ by $\varphi_{i,j}$.
Then the interior angle at the intersection of the two corresponding
planes is clearly $\pi- \varphi_{i,j}$.
Therefore $P(v \in I_{+})$ equals the relative surface area of
a spherical triangle with angles $\pi- \varphi_{1,2}$, $\pi- \varphi
_{1,3}$, $\pi- \varphi_{2,3}$,
%
%
\begin{eqnarray}\label{eqprob3reg}
P(v \in I_{+}) &=& \frac{1}{4\pi} \biggl( \sum
_{1\leq i<j \leq3} (\pi- \varphi_{i,j} ) - \pi\biggr)\nonumber
\\
& =& \frac{1}{4\pi} \biggl(\frac{\pi}{2} + \sum_{1\leq i<j \leq3}
\biggl(\frac{\pi}{2} - \varphi_{i,j} \biggr) \biggr)
\\
&=&
\frac{1}{4\pi} \biggl( \frac{\pi}{2} + \sum
_{1\leq i<j \leq3} \arcsin\biggl( \frac{u_i \cdot u_j}{\llVert u_i
\rrVert \llVert u_j \rrVert } \biggr) \biggr).\nonumber
\end{eqnarray}
By Proposition~\ref{propchtr} we have $c_{i,j} = (\lambda^2 - 3)/6$
and $\varphi_{i,j} = \arccos( (1-\lambda)/4 )$ in the
cherry-transitive case, thus
%
%
\begin{equation}
\label{eqq3} q_3(\lambda) = \frac{1}{8} + \frac{3}{4 \pi}
\arcsin\biggl( \frac
{1-\lambda}{4} \biggr) = \frac{1}{2} -
\frac{3}{4 \pi} \arccos\biggl( \frac{1-\lambda
}{4} \biggr).
\end{equation}
\begin{pf*}{Proof of Theorem~\ref{thm3regmain}}
The statement of the theorem is true for the complete graph $K_4$
as the independence ratio is $1/4$ and the minimum eigenvalue is $-1$
in that case.
For any other $3$-regular transitive graph $G$ we have $\lambda_{\min
}\leq-2$; see Proposition~\ref{proplamin3reg} in the \hyperref[secapp]{Appendix}.
Therefore it suffices to prove that $P(v \in I_{+}) \geq q_3(\lambda
)$, whenever $\lambda\leq-2$.

Recall that $Y_1, Y_2, Y_3$ are standard Gaussians with
pairwise covariances $c_{i,j}$. Therefore the matrix
\[
\pmatrix{ 1 & c_{1,2} & c_{1,3}
\cr
c_{1,2} & 1 &
c_{2,3}
\cr
c_{1,3} & c_{2,3} & 1} %
\]
is positive semidefinite. In particular, its determinant is nonnegative:
\[
1 + 2 c_{1,2} c_{1,3} c_{2,3} - c_{1,2}^2
- c_{1,3}^2 - c_{2,3}^2 \geq0.
\]
Furthermore, according to (\ref{eqcovsum}) we have
$c_{1,2}+c_{1,3}+c_{2,3} = (\lambda^2-3)/2 \geq1/2$, because
$\lambda\leq-2$.
It follows that each $c_{i,j}$ must be between $-1/2$ and $1$.

Indeed, let $x,y,z$ be real numbers between $-1$ and $1$
with $x+y+z \geq1/2$ and $1+2xyz-x^2-y^2-z^2 \geq0$.
Assume that $z < -1/2$. Then
\begin{eqnarray*}
0 &\leq& 1+2xyz-x^2-y^2-z^2 = 1 + 2(z+1)xy -
(x+y)^2 - z^2
\\
& \leq&
1 + 2(z+1) \biggl( \frac{x+y}{2} \biggr)^2 -
(x+y)^2 - z^2 = 1 + \frac{z-1}{2}
(x+y)^2 - z^2
\\
&\leq& 1 + \frac{z-1}{2} \biggl(
\frac{1}{2}-z \biggr)^2 - z^2 < 0,
\end{eqnarray*}
contradiction. Therefore $z \geq-1/2$. Similarly, $x,y \geq-1/2$, too.

Next we bound $u_i \cdot u_j/(\llVert u_i \rrVert \llVert u_j
\rrVert )$ from below.
Using (\ref{eqcovinnerpr}), $x = (y_1+y_2+y_3)/\lambda$
and $c_{1,2}+c_{1,3}+c_{2,3} = (\lambda^2-3)/2$
\begin{eqnarray*}
x \cdot y_1 &=& \frac{1}{\lambda}( 1 + c_{1,2} +
c_{1,3} ) = \frac{1}{\lambda} \biggl( 1 + \frac{\lambda^2-3}{2} -
c_{2,3} \biggr) = \frac{\lambda}{2} - \frac{1}{2\lambda}-
\frac{1}{\lambda} c_{2,3},
\\
\llVert u_1\rrVert^2 &=& \llVert x - y_1
\rrVert^2 = 2 - 2 x \cdot y_1 = 2 - \lambda+
\frac{1}{\lambda} + \frac{2}{\lambda} c_{2,3}.
\end{eqnarray*}
Similar formulas hold for $x \cdot y_i$ and $\llVert u_i\rrVert $, $i=2,3$.
By the inequality of arithmetic and geometric means it follows that
\begin{eqnarray*}
\llVert u_1 \rrVert\llVert u_2 \rrVert&\leq&
\frac{\llVert u_1\rrVert ^2 + \llVert u_2\rrVert ^2}{2} = 2 -
\lambda+ \frac{1}{\lambda} + \frac{1}{\lambda}(
c_{1,3} + c_{2,3} )
\\
&=& \frac{-1}{\lambda} \biggl(
\frac{1}{2} - 2\lambda+ \frac
{\lambda^2}{2} + c_{1,2} \biggr).
\end{eqnarray*}
Note that this holds with equality when all $c_{i,j}$ are equal. Furthermore,
\begin{eqnarray*}
u_1 \cdot u_2 &=& (x-y_1)
\cdot(x-y_2) = 1 + c_{1,2} - x \cdot(y_1+y_2)
\\
&=& 1 + c_{1,2} + x \cdot(y_3 - \lambda x)
\\
&=&
1 + c_{1,2} + \biggl( \frac{\lambda}{2} - \frac{1}{2\lambda
}-
\frac{1}{\lambda} c_{1,2} \biggr) - \lambda
\\
&=& \frac{-1}{\lambda} \biggl(
\frac{1}{2} - \lambda+ \frac{\lambda
^2}{2} + (1-\lambda) c_{1,2}
\biggr).
\end{eqnarray*}
It follows that
\[
\frac{u_1 \cdot u_2}{\llVert u_1\rrVert \llVert u_2\rrVert } \geq
\frac{ (1/2) - \lambda+ (\lambda^2/2) + (1-\lambda
) c_{1,2} }{(1/2) - 2\lambda+ (\lambda^2/2) + c_{1,2} },
\]
because the numerator is positive
(note that $-3 \leq\lambda\leq-2$ and $c_{1,2} \geq-1/2$).
The analogous inequality holds for any other pair of indices $i,j$.
Since $\arcsin$ is a monotone increasing function,
(\ref{eqprob3reg}) yields that
\[
P(v \in I_{+}) \geq\frac{1}{8} + \frac{1}{4\pi} \sum
_{1\leq i < j \leq3} \arcsin\biggl( \frac{ (1/2) - \lambda+ (\lambda^2/2) + (1-\lambda) c_{i,j} }{
(1/2) - 2\lambda+ (\lambda^2/2) + c_{i,j} } \biggr).
\]
Setting
\[
f(t) \stackrel{\mathrm{def}} {=} \arcsin\biggl( \frac{ (1/2) -
\lambda+ (\lambda^2/2) + (1-\lambda) t }{
(1/2) - 2\lambda+ (\lambda^2/2) + t } \biggr),
\]
we have
%
%
\begin{equation}
\label{eqprovedineq} P(v \in I_{+}) \geq\frac{1}{8} +
\frac{1}{4\pi} \sum_{1\leq i < j \leq3} f( c_{i,j} ).
\end{equation}
On the other hand,
%
%
\begin{equation}
\label{eqq3withf} q_3(\lambda) = \frac{1}{8} +
\frac{3}{4\pi} f \biggl( \frac
{\lambda^2-3}{6} \biggr),
\end{equation}
which follows from (\ref{eqq3}) and the definition of $f$.
[It also follows from the fact that when each $c_{i,j}$ is equal to
$(\lambda^2-3)/6$,
then (\ref{eqprovedineq}) should hold with equality.]
In view of (\ref{eqprovedineq}) and (\ref{eqq3withf})
we need to show that
%
%
\begin{equation}
\label{eqineqf} \frac{1}{3} \sum_{1\leq i < j \leq3} f(
c_{i,j} ) \geq f \biggl( \frac{\lambda^2-3}{6} \biggr),
\end{equation}
where each $c_{i,j}$ is between $-1/2$ and $1$, and their average is
$(\lambda^2-3)/6$.
This, of course, would follow from the convexity of $f$.
Unfortunately, $f$ is not convex on the entire interval $[-1/2,1]$.
We claim, however, that the tangent line to $f$ at $t_0 = (\lambda^2-3)/6$
is below $f$ on the entire interval $[-1/2,1]$, which still implies
(\ref{eqineqf}).
The rather technical proof of this claim can be found
in the \hyperref[secapp]{Appendix} (Lemma~\ref{lemtangentline}).

Now let $\lambda= \lambda_{\min}\leq-2$, then $P(v \in I_{+}) \geq
q_3(\lambda_{\min})$.
So the expected size of the random independent set $I_{+}$ is at least
$q_3(\lambda_{\min}) \llvert V(G)\rrvert $;
thus the independence ratio of $G$ is at least $q_3(\lambda_{\min})$.

To prove the second part of the statement we notice that
the random independent set $I_{-}$ (see Definition~\ref{defindset})
has the same expected size. Indeed, if we replace $X_v$, $v \in V(G)$ with
$X'_v = - X_v$, then $X'_v$, $v \in V(G)$ have the same joint distribution
and the roles of $I_{+}$ and $I_{-}$ interchange.
Since $I_{+}$ and $I_{-}$ are always disjoint,
the expected size of their union $I_{+} \cup I_{-}$ is at least $2
q_3(\lambda_{\min}) \llvert V(G)\rrvert $.
Consequently, there must exist disjoint independent sets $I_1,I_2$ in $G$
with $\llvert I_1 \cup I_2\rrvert \geq2 q_3(\lambda_{\min}) \llvert
V(G)\rrvert $.
\end{pf*}

For graphs with very large odd-girth, Theorem~\ref{thm3regother}
of the \hyperref[secapp]{Appendix} gives a slightly better bound.
The proof is based on the same random eigenvector,
but uses a different method to find large independent sets.

\subsection{The arc-transitive case} \label{subsecedgetr}

The following innocent-looking, and very plausible, conjecture
is open in dimension $n \geq4$.
%

\begin{conjecture} \label{conjgeom}
Let $S$ be a sphere in the $n$-dimensional Euclidean space~$\mathbb{R}^n$.
We have $n+1$ spherical caps with the same given radius on $S$.
We want to find the configuration for which
the volume of the union of the caps is maximal.
It is conjectured that this optimal configuration is
always the one where the $n+1$ centers are
the vertices of a regular simplex in $\mathbb{R}^n$.
\end{conjecture}

The statement of the conjecture is trivial for $n=2$,
while the $n=3$ case follows from the so-called \textit{moment theorem}
of Fejes T\'oth \cite{fejestoth}, Theorem~2; see also \cite{regularfigures}, Section~34.
The genereal case would follow from the following conjecture:
the volume of the intersection of a fixed spherical cap and a spherical simplex
of fixed volume is maximal when the spherical simplex is regular, and
its center coincides with the center of the spherical cap
(G\'abor Fejes T\'oth, personal communication,~2012).

In what follows we will explain
how the case $n=d-1$ of Conjecture~\ref{conjgeom}
implies that $P( v \in I_{+} ) \geq q_d(\lambda)$ holds
for every $d$-regular arc-transitive graph $G$,
and consequently the independence ratio of $G$ is at least $q_d(\lambda
_{\min})$.
In particular, the $d=4$ case follows from the $n=3$ case of
the conjecture, which is known to be true; see Theorem~\ref{thmarctr}.
Using our previous notation, $P( v \in I_{+} )$ is the volume
of the spherical simplex $T$ determined by the half-spaces
with outer normal vectors $-u_i$, $i=1, \ldots, d$,
while $q_d(\lambda)$ is the volume of the same simplex
in the case when all the angles
$\varphi_{i,j} = \angle(u_i,u_j)$, $i \neq j$
are the same.
In other words, we need to show that
the volume of the spherical simplex $T$
is minimal when the angles $\angle(u_i,u_j)$ are the same.

If $G$ is arc-transitive, then the covariances
$\cov(X,Y_i) = x \cdot y_i$ are all equal. Since
\[
x \cdot y_1 + \cdots+ x \cdot y_d = x
\cdot(y_1+ \cdots+ y_d) = x \cdot(\lambda x) = \lambda,
\]
we get that $x \cdot y_i = \lambda/d$ for each $i$.
It follows that the angle enclosed by $x$ and $u_i$
%
%
\begin{equation}
\label{eqdelta} \angle(x,u_i) = \delta\stackrel{\mathrm{def}} {=}
\frac{ \pi-
\arccos( \lambda/ d)
}{2} = \frac{ \arccos( - \lambda/d ) }{2}\qquad\mbox{for each } i.
\end{equation}

Now let $S_l$ be the set of points on $S^{d-1}$ that
has some fixed distance $l$ from $x$;
thus $S_l$ is a $(d-2)$-dimensional sphere for any $l$.
The intersection of $S_l$ and the half-space
with outer normal vector $u_i$ is
a spherical cap of radius depending only on $l$ and $\lambda$.
So the intersection of $S_l$ and our spherical simplex $T$ can be obtained
by removing $d$ spherical caps of the same given radius from $S_l$.
If Conjecture~\ref{conjgeom} is true for $n=d-1$,
then the total volume of the removed area is maximal
for the ``regular configuration'' when each $\angle(u_i,u_j)$ is the same.
Therefore the $(d-2)$-dimensional volume of $T \cap S_l$
is minimal for the regular configuration for any $l$.
It follows that the $(d-1)$-dimensional volume of $T$
is also minimal for the regular configuration,
and this is what we wanted to prove.

\subsection{Bounds near $-d$} \label{secnearnegd}

Even if Conjecture~\ref{conjgeom} is not assumed to be true,
the above observations yield a lower bound
for the independence ratio of $d$-regular arc-transitive graphs
in the case when the least eigenvalue is close to $-d$.

\begin{pf*}{Proof of Theorem~\ref{thmarctr}}
We already explained in Section~\ref{subsecedgetr} why
Conjecture~\ref{conjgeom} implies that the independence ratio is at
least $q_d(\lambda_{\min})$.

Next we prove tht first part of Theorem~\ref{thmarctr}.
As we have seen in (\ref{eqdelta}), $ \angle(x,u_i) = \delta$ for
each $i$,
which means that each point of $S^{d-1}$ at (spherical) distance
less than $\pi/2 - \delta$ from $x$ is contained in our spherical
simplex $T$.
These points form a spherical cap with center $x$ and radius $\pi/2 -
\delta$.
(In fact, this spherical cap is the ``inscribed ball'' of $T$.)
Using (\ref{eqdelta}) and that
$\arccos(t) \leq\pi/2 \sqrt{1-t}$ for any $t \in[0,1]$, we get
\[
\delta= \frac{ \arccos( - \lambda/d ) }{2} \leq\frac{\pi}{4} \sqrt
{1+\lambda/d} %
\]
provided that $\lambda\leq0$.

This spherical cap can be obtained by taking the hemisphere (around $x$),
and removing a strip of ``width'' $\delta$ (in spherical distance).
The volume of this strip is clearly at most $\delta\vol(S^{d-2})$. Therefore
the volume of the spherical cap is at least $ \vol(S^{d-1})/2 -
\delta\vol(S^{d-2})$, whence
\[
P( v \in I_{+} ) \geq\frac{ \vol(S^{d-1})/2 - \delta\vol(S^{d-2})
}{ \vol(S^{d-1}) } = \frac{1}{2} -
\frac{ \pi\vol(S^{d-2}) }{ 4 \vol(S^{d-1}) } \sqrt{\frac{\lambda
+d}{d}}.
\]
For $d=4$ we have $\vol(S^{2}) / \vol(S^{3}) = (4 \pi) / (2 \pi^2)
= 2/\pi$,
so the bound is
\[
\tfrac{1}{2} - \tfrac{1}{4} \sqrt{\lambda+4}.
\]
For general $d$, we use the estimate
$\vol(S^{d-2}) / \vol(S^{d-1}) \leq\sqrt{d}/\sqrt{2\pi}$
(see Lemma~\ref{lemvolratio} of the \hyperref[secapp]{Appendix}) to obtain the
following bound:
\[
\frac{1}{2} - \frac{\sqrt{\pi}}{4\sqrt{2}} \sqrt{\lambda+d} > \frac
{1}{2} -
\frac{1}{3} \sqrt{\lambda+d}.
\]
These are lower bounds for the probability $P( v \in I_{+} )$,
in particular, for $q_d(\lambda)$.
Thus the first part of Theorem~\ref{thmarctr} follows,
as well as the estimate (\ref{eqq4}) for $q_4(\lambda)$
and the last statement of Proposition~\ref{propqdprop}.
\end{pf*}

\begin{pf*}{Proof of Theorem~\ref{thmcrude}}
In the general (vertex-transitive) case, we use that
$ x \cdot y_1 + \cdots+ x \cdot y_d = \lambda$ and $x \cdot y_j \geq-1$:
\[
x \cdot y_i \leq\lambda+ d - 1\qquad\mbox{for each } 1 \leq i \leq d.
\]
Therefore the angle $\angle(x,y_i)$ is at least $\arccos( \lambda+
d - 1 )$.
Using that $\arccos(t) \leq\pi/2 \sqrt{1-t}$ for any $t \in[0,1]$,
it follows that
\[
\angle(x,u_i) \leq\delta' \stackrel{\mathrm{def}} {=}
\frac{ \pi-
\arccos( \lambda
+ d - 1 ) }{2} = \frac{ \arccos( 1 - \lambda- d ) }{2} \leq\frac{\pi
}{4} \sqrt{\lambda+d}
\]
provided that $\lambda\leq-d+1$.
This means that our spherical simplex $T$ contains
the spherical cap with center $x$ and radius $\pi/2 - \delta'$. Therefore
\begin{eqnarray*}
P( v \in I_{+} ) &\geq&\frac{ \vol(S^{d-1})/2 - \delta' \vol
(S^{d-2}) }{ \vol(S^{d-1}) } = \frac{1}{2} -
\frac{ \pi\vol(S^{d-2}) }{ 4 \vol(S^{d-1}) } \sqrt{\lambda+d}
\\
&\geq& \frac
{1}{2} - \frac{\sqrt{\pi}}{4\sqrt{2}}
\sqrt{d(\lambda+d)}.
\end{eqnarray*}
Since $\sqrt{\pi}/(4\sqrt{2}) < 1/3$, Theorem~\ref{thmcrude} follows.
\end{pf*}


\section{Infinite transitive graphs} \label{sec3}

\subsection{Random wave functions}

Our goal now is to generalize the random eigenvectors
we introduced in Section~\ref{sec2} for infinite transitive graphs $G$.
For an infinite graph $G$ the adjacency operator
$A_G\dvtx  \ell_2( V(G) ) \to\ell_2( V(G) )$ might not have any eigenvectors
(i.e., the point spectrum might be empty).
So the approach we used in the finite setting will not work here.
Instead, we will define random wave functions as the limit of linear
factor of i.i.d. processes.
The coefficients of these linear factors will be approximate
eigenvectors of $A_G$
that are invariant under automorphisms fixing some root $x \in V(G)$.
We start with proving that such approximate eigenvectors exist
for any $\lambda$ in the spectrum $\sigma(A_G)$.
Let $\St_x(G)$ denote the \emph{stabilizer subgroup},
that is, the group of automorphisms fixing $x$.
%

\begin{theorem} \label{thminvapprev}
Let $G$ be an infinite vertex-transitive graph
with adjacency operator $A_G$ and with some fixed root $x$.
Then for any $\varepsilon>0$ and any $\lambda$ in the spectrum
$\sigma(A_G)$, there exists
a $\St_x(G)$-invariant vector $\alpha\in\ell_2( V(G) )$ such that
\[
\llVert\alpha\rrVert=1\quad\mbox{and}\quad \llVert A_G \alpha- \lambda
\alpha\rrVert\leq\varepsilon.
\]
\end{theorem}

\begin{pf}
Consider the projection-valued measure $P_\lambda$
corresponding to the self-adjoint operator $A_G$.
This ``measure'' assigns an orthogonal projection
$P_S$ to each Borel set $S \subseteq\mathbb{R}$.
According to spectral theory, one can integrate with respect to this measure.
For instance, the following formula holds:
\[
A_G = \int_\mathbb{R}\lambda\,\mathrm{d}
P_\lambda.
\]

Furthermore, the projections $P_S$ have the property that if an
operator $T$
commutes with $A_G$, then it also commutes with each projection $P_S$.
There is a unitary operator $U_\Phi$ corresponding to each $\Phi\in
\Aut(G)$
[the one that permutes the coordinates of $\ell_2( V(G) )$ according
to $\Phi$].
Since $U_\Phi$ commutes with $A_G$, it also commutes with the
projections $P_S$.

Now let $\lambda_0$ be an arbitrary element of the spectrum $\sigma(A_G)$,
and set $S=[\lambda_0-\varepsilon, \lambda_0+\varepsilon]$. We
define $\alpha$ as the image of
the indicator function $\ind_x$ under the projection $P_S$,
\[
\alpha\stackrel{\mathrm{def}} {=}P_{[\lambda_0-\varepsilon, \lambda
_0+\varepsilon]} \ind_x.
\]
Note that $\ind_x$ is a fixed point of $U_\Phi$ for any $\Phi\in\St
_x(G)$. Therefore
\[
U_\Phi\alpha= U_\Phi P_{[\lambda_0-\varepsilon, \lambda
_0+\varepsilon]} \ind_x
= P_{[\lambda_0-\varepsilon, \lambda_0+\varepsilon]} U_\Phi\ind_x =
P_{[\lambda_0-\varepsilon, \lambda_0+\varepsilon]}
\ind_x = \alpha, %
\]
and thus $\alpha$ is $\St_x(G)$-invariant.
On the other hand, since $P_S P_{\mathbb{R}\setminus S} = 0$, we have
\[
A_G \alpha- \lambda_0 \alpha= \biggl( \int
_{\mathbb{R}} (\lambda-\lambda_0) \,\mathrm{d}
P_\lambda\biggr) \alpha= \biggl( \int_{[\lambda_0-\varepsilon, \lambda
_0+\varepsilon]} (\lambda-
\lambda_0) \,\mathrm{d} P_\lambda\biggr) \alpha,
\]
which clearly implies that
\[
\llVert A_G \alpha- \lambda_0 \alpha\rrVert\leq
\varepsilon\llVert\alpha\rrVert.
\]
It remains to show that $\alpha= P_{[\lambda_0-\varepsilon,
\lambda_0+\varepsilon]} \ind_x \neq0 $.
Assume that $P_{[\lambda_0-\varepsilon, \lambda_0+\varepsilon]}
\ind_x = 0$.
It follows that $P_{[\lambda_0-\varepsilon, \lambda
_0+\varepsilon]} \ind_v = 0$ for every vertex $v \in V(G)$.
Indeed, let $\Phi\in\Aut(G)$ such that $\Phi x = v$. Then $U_\Phi
\ind_x = \ind_v$ and
\[
P_{[\lambda_0-\varepsilon, \lambda_0+\varepsilon]} \ind_v = P_{[\lambda
_0-\varepsilon, \lambda_0+\varepsilon]} U_\Phi\ind
_x = U_\Phi P_{[\lambda_0-\varepsilon, \lambda_0+\varepsilon]} \ind_x
= 0.
\]
This holds for each vertex $v$, which clearly implies
that $P_{[\lambda_0-\varepsilon, \lambda_0+\varepsilon]} = 0$.
Then the operator
\[
B = \int_{\mathbb{R}\setminus[\lambda_0-\varepsilon, \lambda
_0+\varepsilon]} \frac{1}{\lambda-\lambda_0} \,\mathrm{d} P_\lambda
\]
would be the inverse of $A_G - \lambda_0 I$ contradicting
our assumption that $\lambda_0 \in\sigma(A_G)$.
\end{pf}

%
\begin{remark}
There is a general theorem for Hilbert spaces
saying that every point of the spectrum of a self-adjoint operator
is an approximate eigenvalue \cite{sunder}, Corollary 4.1.3.
So the real content of the above theorem is that one can find
approximate eigenvectors that are $\St_x(G)$-invariant.
This invariance will be crucial for us later on,
when we will use these approximate eigenvectors
as coefficients to define linear factor of i.i.d. processes.
\end{remark}

Suppose now that we have an i.i.d. process on $G$:
independent standard normal random variables $Z_u$ assigned to each
vertex $u$.
We will consider processes $X_v$, $v \in V(G)$, where each $X_v$ is
a (possibly infinite) linear combination of $Z_u$, $u \in V(G)$.
We collected some obvious properties of such processes in the next proposition.
%

\begin{proposition} \label{proplinfactor}
Let $\beta_{v,u}$, $v,u \in V(G)$ be real numbers, and let
%
%
\begin{equation}
\label{eqlin} X_v = \sum_{u \in V(G)}
\beta_{v,u} Z_u.
\end{equation}
The infinite sum in (\ref{eqlin}) converges almost surely if and only if
%
%
\begin{equation}
\label{eqsquaresum} \sum_{u \in V(G)} \beta_{v,u}^2
< \infty.
\end{equation}
If (\ref{eqsquaresum}) is satisfied,
then $X_v$ is a centered Gaussian with variance
$\var(X_v) = \sum_{u \in V(G)} \beta_{v,u}^2$.

The process $X_v$, $v \in V(G)$ is $\Aut(G)$-invariant if and only if
%
%
\begin{equation}
\label{eqinv} \beta_{v,u} = \beta_{\Phi v, \Phi u}\qquad\mbox{for all }
\Phi
\in\Aut(G).
\end{equation}
\end{proposition}

Now we are in a position to formally define linear factor of i.i.d. processes.
%

\begin{definition} \label{deflinfactor}
We say that a process $X_v$, $v \in V(G)$ is a
\emph{linear factor} of the i.i.d. process $Z_u$ if
it can be written as in (\ref{eqlin})
for some real numbers $\beta_{v,u}$, $v,u \in V(G)$
satisfying (\ref{eqsquaresum}) and (\ref{eqinv}).
\end{definition}
%

\begin{remark} \label{rmlinfactor}
Let us fix a root $x \in V(G)$. For a linear factor the coefficients
$\alpha_u \stackrel{\mathrm{def}}{=}\beta_{x,u}$ clearly determine
each $\beta_{v,u}$.
Here $\alpha=(\alpha_u)_{u \in V(G)}$ can be any $\St
_x(G)$-invariant vector in $\ell_2( V(G) )$.
So there is a one-to-one correspondance between linear factor of
i.i.d. processes on $G$
and $\St_x(G)$-invariant vectors $\alpha\in\ell_2( V(G) )$.
Also, by Proposition~\ref{proplinfactor} we have $\var(X_v) = \llVert
\alpha\rrVert ^2$.
\end{remark}

Recall Definition~\ref{defgaussianprocess} of invariant Gaussian processes.
%

\begin{definition} \label{defgaussianev}
We call an invariant Gaussian process $X_v$, $v \in V(G)$
a~\emph{Gaussian wave function} with eigenvalue $\lambda$ if
\[
\sum_{u \in N(v)} X_u = \lambda
X_v\qquad\mbox{for each vertex } v \in V(G),
\]
where $N(v)$ denotes the set of neighbors of $v$ in $G$.
\end{definition}

%
\begin{example}
It was shown in \cite{csghv} that for the $d$-regular tree $T_d$
there exists an essentially unique Gaussian wave function for each
$\lambda\in[-d,d]$.
Furthermore, this Gaussian wave function can be approximated
by factor of i.i.d. processes provided that
$\lambda$ is in the spectrum $\sigma(T_d) = [-2\sqrt{d-1}, 2\sqrt{d-1}]$.
\end{example}

In general, it is not clear for which $\lambda$ such Gaussian wave
functions exist
and whether they are unique.
%

\begin{definition}
For a transitive graph $G$ we call the closed set
\[
\widetilde{\sigma}(G) \stackrel{\mathrm{def}} {=} \{\mbox{$\lambda$: there
exists a Gaussian wave function on $G$ with eigenvalue $\lambda$} \} %
\]
the \emph{Gaussian spectrum} of $G$.
\end{definition}

Theorem~\ref{thmgaussianev} claims that
for any $\lambda\in\sigma(A_G)$ there exists a Gaussian wave
function on $G$,
which can be approximated by linear factor of i.i.d. processes.
Therefore $\widetilde{\sigma}(G) \supseteq\sigma(A_G)$.

\begin{pf*}{Proof of Theorem~\ref{thmgaussianev}}
We use the $\St_x(G)$-invariant approximate eigenvectors
of Theorem~\ref{thminvapprev} to define linear factor of i.i.d. processes.
So let $\varepsilon> 0$ be arbitrary and $\alpha^\varepsilon$ a
$\St_x(G)$-invariant vector with
$\llVert \alpha^\varepsilon\rrVert = 1$ and $\llVert A_G \alpha
^\varepsilon-
\lambda\alpha^\varepsilon\rrVert \leq\varepsilon$.
By Remark~\ref{rmlinfactor} for each $\alpha^\varepsilon$ there is
a corresponding linear factor $X_v^\varepsilon$, $v \in V(G)$. Note
that the process
\[
Y_v^\varepsilon\stackrel{\mathrm{def}} {=}\sum
_{u \in N(v)} X_u^\varepsilon- \lambda
X_v^\varepsilon%
\]
is also a linear factor, and the corresponding
coefficient vector is $\delta^\varepsilon\stackrel{\mathrm
{def}}{=}A_G \alpha
^\varepsilon- \lambda\alpha^\varepsilon$.
Therefore $X_v^\varepsilon$ is an invariant Gaussian process with
$\var(X_v^\varepsilon) = \llVert \alpha^\varepsilon\rrVert ^2 = 1$ and
\[
\var\biggl( \sum_{u \in N(v)} X_u^\varepsilon-
\lambda X_v^\varepsilon\biggr) = \var\bigl(
Y_v^\varepsilon\bigr) = \bigl\llVert\delta^\varepsilon\bigr
\rrVert^2 = \bigl\llVert A_G \alpha^\varepsilon-
\lambda\alpha^\varepsilon\bigr\rrVert^2 \leq\varepsilon^2.
\]
Since the space of invariant Gaussian processes with variance $1$ is compact,
it follows that there exists a sequence $\varepsilon_n$ converging to $0$
such that the processes $X_v^{\varepsilon_n}$ converge in distribution.
The limit process will be a nontrivial invariant Gaussian process $X_v$
that satisfies the eigenvector equation (\ref{eqeigen}) at each vertex.
\end{pf*}

\subsection{Factor of i.i.d. processes}

For a graph $G$ we defined an i.i.d. process on $G$
as independent standard normal random variables $Z_v$, $v \in V(G)$.
In other words, $Z = ( Z_v )_{v \in V(G)}$ is
a random point in the measure space $(\Omega, \mu)$,
where $\Omega$ is $\mathbb{R}^{V(G)}$ with the product topology,
and $\mu$ is the product of standard Gaussian measures (one on each
copy of $\mathbb{R}$).
The natural action of $\Aut(G)$ on $V(G)$ gives rise to
an action of $\Aut(G)$ on $\Omega$: for $\Phi\in\Aut(G)$
and $\omega= ( \omega_v )_{v \in V(G)} \in\Omega$ let
\[
( \Phi\cdot\omega)_v \stackrel{\mathrm{def}} {=}\omega
_{ \Phi^{-1} v }.
\]
Let $G$ be an infinite transitive graph, and suppose that
$F$ is a measurable $\Omega\to\Omega$ function
that is $\Aut(G)$-equivariant [i.e., commutes with the $\Aut(G)$-action].
Then $X = F(Z)$ is an invariant process on $G$.
Such a process $X = ( X_v )_{v \in V(G)}$ is called
a \emph{factor} of the i.i.d. process $Z$.

An $\Aut(G)$-equivariant $F \colon\Omega\to\Omega$ function
is determined by $f = \pi_x \circ F$, where $\pi_x \colon\Omega\to
\mathbb{R}$ is
the projection corresponding to the coordinate of some fixed root $x$.
Here $f$ can be any $\St_x(G)$-invariant $\Omega\to\mathbb{R}$ function.
So factor of i.i.d. processes can be identified with
measurable, $\St_x(G)$-invariant functions $f \colon\Omega\to
\mathbb{R}$.

Next we will prove Theorems~\ref{thmnotclosed}~and~\ref{thmcayley}
by showing that certain Gaussian wave functions $X_v$, $v \in V(G)$
cannot be obtained as factor of i.i.d. processes.
Since $X_v$ has finite variance in that case,
we can restrict ourselves to functions $f \in L_2(\Omega, \mu)$.
Let $H_{\mathrm{inv}}\subset L_2(\Omega, \mu)$ be the subspace containing
those $f \in L_2(\Omega, \mu)$ that are $\St_x(G)$-invariant.
There is a natural way to define an adjacency operator~$\mathcal{A}$
on the Hilbert space $H_{\mathrm{inv}}$. Let
\[
( \mathcal{A}f ) (\omega) \stackrel{\mathrm{def}} {=} \sum
_{y \in N(x)} f ( \Phi_{y \to x} \cdot\omega),
\]
where $\Phi_{y \to x}$ is an (arbitrary) automorphism of $G$ taking
$y$ to $x$.
Since $f$ is $\St_x(G)$-invariant, $\mathcal{A}$ is well defined.

Suppose now that we have a Gaussian wave function with eigenvalue
$\lambda$
that can be obtained as a factor of i.i.d. process.
Then the corresponding $f$ satisfies the eigenvector equation $\mathcal
{A}f = \lambda f$.
In particular, $\lambda$ needs to be in the point spectrum of
$\mathcal{A}$.
[Note that an eigenvector $f$ of $\mathcal{A}$ does not necessarily
give us a
Gaussian wave function: although the corresponding factor of i.i.d. process
will satisfy the eigenvector equation at each vertex,
$f(Z)$ might not have a Gaussian distribution.]

\begin{pf*}{Proof of Theorem~\ref{thmnotclosed}}
Since $L_2(\Omega, \mu)$ is a separable Hilbert space,
so is $H_{\mathrm{inv}}$, and consequently the point spectrum of
$\mathcal{A}\dvtx  H_{\mathrm{inv}}
\to H_{\mathrm{inv}}$ is countable.

Therefore only for countably many $\lambda$'s
can we have a Gaussian wave function on $G$ that
can be obtained as a factor of i.i.d. process.
However, if $\sigma(A_G)$ is uncountable, then by Theorem~\ref{thmgaussianev},
$G$ has Gaussian wave functions for uncountably many different
eigenvalues $\lambda$;
moreover, they can all be approximated by linear factor of i.i.d. processes.
\end{pf*}

\begin{pf*}{Proof of Theorem~\ref{thmcayley}}
We will use two basic facts about the point spectra of
the adjacency operators $A_G$ and $\mathcal{A}$.
First, $\lambda_{\max}$ is never in the point spectrum $\sigma_p(A_G)$
(we will give a short proof for this in the \hyperref[secapp]{Appendix}; see Lemma~\ref
{lemlama}).
Second, $\sigma_p(\mathcal{A}) \subseteq\sigma_p(A_G) \cup\{d\}$
for Cayley graphs
(this will be explained after the proof).
Therefore $\lambda_{\max}$ is not in the point spectrum of $\mathcal{A}$
provided that \mbox{$\lambda_{\max}<d$}, and consequently,
a Gaussian wave function with eigenvalue $\lambda_{\max}$
cannot be obtained as a factor of i.i.d. process.

In the case $\lambda_{\max}= d$
the Gaussian wave function has to be constant;
that is, $X_u = X_v$ for any two vertices $u,v$.
However, for a factor of i.i.d. process the correlation between $X_u$
and $X_v$
should tend to $0$ as the distance of $u$ and $v$ goes to infinity; see
Proposition~\ref{propcorrdecay} in the \hyperref[secapp]{Appendix}.
\end{pf*}

Next we will explain the relation between
the adjacency operators $A_G$ and $\mathcal{A}$.
This can be found in \cite{kechristsankov}, Section~3, in a more
general setting;
see also \cite{lyonsnazarov}, Theorem~2.1 and Corollary 2.2.
Let $\nu$ denote the standard Gaussian measure.
Since $L_2(\mathbb{R}, \nu)$ is a separable Hilbert space,
it has a countable orthonormal basis: $g_0, g_1, g_2, \ldots,$
where $g_0$ will be assumed to be the constant $1$ function.
Let $\mathcal{I}$ denote the set of finitely supported $V(G) \to\{
0,1,2,\ldots\}$ functions.
For each $q \in\mathcal{I}$ we define an $\Omega\to\mathbb{R}$ function:
\[
W_q(\omega) \stackrel{\mathrm{def}} {=}\prod
_{v \in V(G)} g_{q(v)} ( \omega_v ).
\]
Note that this is actually a finite product,
since all but finitely many terms are equal to $g_0 \equiv1$.
According to \cite{kechristsankov}, Lemma~3.1,
the functions $W_q$, $q\in\mathcal{I}$ form an orthonormal basis of
$L_2(\Omega, \mu)$.
It follows that $L_2(\Omega, \mu)$ is separable,
which fact was used in the proof of Theorem~\ref{thmnotclosed}.

We defined the operator $\mathcal{A}$ on the space $H_{\mathrm
{inv}}\subset
L_2(\Omega, \mu)$
containing\break $\St_x(G)$-invariant functions.
When $G$ is a Cayley graph, there is a natural way to
extend $\mathcal{A}$ to an adjacency operator over the whole space
$L_2(\Omega, \mu)$.
Suppose that $\Gamma$ is a finitely generated infinite group.
Let $S$ be a finite, symmetric set of generators,
and let $G$ be the corresponding Cayley graph;
that is, $V(G) = \Gamma$ and the vertex $v \in\Gamma$
is adjacent to the vertices $\gamma v$, $\gamma\in S$.
The natural action of $\Gamma$ on itself gives rise to
the following $\Gamma$-action on $\Omega$:
\[
( \gamma\cdot\omega)_v \stackrel{\mathrm{def}} {=} ( \omega
)_{\gamma^{-1} v}.
\]
(This is often called the \emph{generalized Bernoulli shift}.)
Then for $f \in L_2(\Omega, \mu)$, let
\[
( \mathcal{A}f ) (\omega) \stackrel{\mathrm{def}} {=} \sum
_{\gamma\in S} f ( \gamma\cdot\omega).
\]
This clearly extends our earlier definition of $\mathcal{A}$.

There is a natural $\Gamma$-action on $\mathcal{I}$ as well: for $q
\in\mathcal{I}$
\[
( \gamma\cdot q ) (v) \stackrel{\mathrm{def}} {=}q \bigl( \gamma^{-1}
v \bigr).
\]
It is compatible with the $\Gamma$-action on $\Omega$ in the
following sense:
\[
W_{\gamma\,\cdot\, q}(\omega) = W_q \bigl( \gamma^{-1} \cdot
\omega\bigr).
\]
It means that
%
%
\begin{equation}
\label{eqAop} \mathcal{A}W_q = \sum_{\gamma\in S}
W_{\gamma\,\cdot\, q}.
\end{equation}
We now consider the orbit $\{ \gamma\cdot p\dvtx  \gamma\in
\Gamma\}$
of a given element $p \in\mathcal{I}$ and
the closure of the space spanned by the corresponding functions
$W_{\gamma\cdot p}$,
\[
H_p \stackrel{\mathrm{def}} {=}\cl\bigl( \spn\{ W_{\gamma
\cdot p}\dvtx  \gamma\in\Gamma\} \bigr) \subset L_2(\Omega, \mu).
\]
It is clear from (\ref{eqAop}) that $H_p$ is $\mathcal{A}$-invariant.
If $p \equiv0$, then $H_p$ consists of the constant functions on
$\Omega$
and both the point spectrum and the spectrum of $ \mathcal{A}\mid
_{H_p} $ is~$\{ d \}$.
Otherwise the stabilizer $\Gamma_p$ of $p$ is a\vspace*{1pt} finite subgroup of
$\Gamma$, and
$ \mathcal{A}\mid_{H_p} $ is closely related to the
original adjacency operator $A_G$.
Indeed, let $T_p \colon H_p \to\ell_2( V(G) ) \cong\ell_2(\Gamma
)$ be the operator defined by
\[
T_p \colon W_q \mapsto\ind_{ \{\gamma\in\Gamma\dvtx  \gamma
\cdot p = q\} },
\]
where $q$ is in the orbit of $p$.
It is easy to see that $T_p$ is a bounded operator
for which $T_p \mathcal{A}\mid_{H_p} = A_G T_p$.
Since $T_p$ is also bounded below, it follows that
\[
\sigma( \mathcal{A}\mid_{H_p} ) \subseteq\sigma(A_G)
\quad\mbox{and}\quad \sigma_p ( \mathcal{A}\mid_{H_p} ) \subseteq
\sigma_p(A_G) %
\]
with equality when the stabilizer $\Gamma_p$ is trivial.

Therefore for Cayley graphs the operators
$A_G \colon\ell_2( V(G) ) \to\ell_2( V(G) )$ and
$\mathcal{A}\colon L_2(\Omega, \mu) \to L_2(\Omega, \mu)$
have the same spectra and point spectra
with the possible exception of the point $d$,
\[
\sigma(\mathcal{A}) = \sigma(A_G) \cup\{d\}\quad\mbox{and}\quad
\sigma_p(\mathcal{A}) = \sigma_p(A_G) \cup
\{d\}.
\]
Consequently,
\[
\sigma_p ( \mathcal{A}\mid_{H_{\mathrm{inv}}} ) \subseteq
\sigma_p(\mathcal{A}) = \sigma_p(A_G) \cup
\{d\},
\]
which we used in the proof of Theorem~\ref{thmcayley}.

\subsection{Independent sets}

Let $G$ be an infinite transitive graph
and $\lambda_{\min}$ be the minimum of its spectrum $\sigma(A_G)$.
Consider linear factor of i.i.d. processes $X_v^n$
converging in distribution to a Gaussian wave function $X_v$
with eigenvalue $\lambda_{\min}$ as $n \to\infty$ as in Theorem
\ref{thmgaussianev}.
We define the following independent sets on $G$:
\[
I_{+} \stackrel{\mathrm{def}} {=} \bigl\{ v\dvtx  X_v >
X_u, \forall u \in N(v) \bigr\}\quad\mbox{and}\quad I_{+}^n
\stackrel{\mathrm{def}} {=} \bigl\{ v\dvtx  X_v^n >
X_u^n, \forall u \in N(v) \bigr\}.
\]
Then for each $n$ the independent set $I_{+}^n$ is a factor of
the i.i.d. process $Z_v$; that is, it is obtained as a measurable function
of $Z_v$, $v \in V(G)$ that commutes with the natural action of $\Aut(G)$.
Furthermore, since the event $v \in I_+$ corresponds to an open set, we have
\[
\liminf_{n \to\infty} P\bigl( v \in I_+^n \bigr) \geq P(
v \in I_+ ).
\]
Therefore whenever we have a lower bound $q$ for $P( v \in I_+ )$,
it yields that for any $\varepsilon>0$ there exists a factor of
i.i.d. independent set with ``size'' greater than $q - \varepsilon$.

Bounding $P( v \in I_+ )$, however, leads us to the same
optimization problem as in the finite case.
We need to estimate the volume of the same spherical simplex
with the exact same constraints.
[Of course, there might be a difference between the finite and infinite setting
in terms of what covariances $c_{i,j}$ can actually come up,
but our proofs used only the trivial constraints
that they form a positive semidefinite matrix and their sum is
$(\lambda_{\min}^2-d)/2$,
which are true in the infinite case, too.]
Thus we obtain the exact same bounds, and Theorem~\ref{thminf} follows.

Actually, in Theorem~\ref{thm3regmain} we proved the bound
only for graphs with $\lambda_{\min}\leq-2$ and argued that
the only finite, $3$-regular, transitive graph for which this does not hold
is the complete graph $K_4$. For infinite transitive graphs
$\lambda_{\min}\leq-2$ holds with no exception. This follows from
the fact
that they contain arbitrarily long paths as induced subgraphs.


\begin{appendix}
\section*{Appendix} \label{secapp}

\setcounter{theorem}{0}

%
%
\begin{theorem} \label{thm3regother}
Suppose that $G$ is a finite, $3$-regular, vertex-transitive graph
with minimum eigenvalue $\lambda_{\min}$ and odd-girth $g$. Then the
independence ratio of $G$ is at least
\[
\frac{5g-3}{16g} + \frac{g+1}{2g} \frac{3}{4 \pi} \arcsin\biggl(
\frac{\lambda_{\min}^2-3}{6} \biggr) \geq\frac{5}{16} + \frac{3}{8 \pi
} \arcsin
\biggl( \frac{\lambda_{\min
}^2-3}{6} \biggr) - \frac{3}{16g}.
\]
In fact, there exist two disjoint independent sets in $G$ such that
their average size divided by $\llvert V(G)\rrvert $ is not less than
the above bound.
\end{theorem}

\begin{pf} 
It is easy to check the statement for $K_4$.
According to Proposition~\ref{proplamin3reg}
$\lambda_{\min}\leq-2$ holds for any other finite, $3$-regular,
transitive graph $G$.
Let $X_v$, $v \in V(G)$ be the random eigenvector corresponding to
$\lambda_{\min}$.
Let $V_{+}$ denote the set of ``positive vertices,'' that is,
\[
V_{+} \stackrel{\mathrm{def}} {=} \bigl\{ v\in V(G)\dvtx
X_v > 0 \bigr\}.
\]
The expected size of $V_{+}$ is $\llvert V(G)\rrvert /2$.

Since $\lambda_{\min}$ is negative, a vertex and its three neighbors
cannot all be positive.
Therefore each vertex has degree at most two in the induced subgraph $G[V_{+}]$.
Thus each connected component of this subgraph is a path or a cycle.
We want to choose an independent set from each component.
We can choose at least half the vertices from paths and even cycles.
From an odd cycle of length $l \geq g$ we can choose $(l-1)/2$ vertices,
which is at least a $(g-1)/(2g)$ proportion of all vertices in that component.
(Recall that $g$ denotes the odd-girth of $G$, i.e.,
the length of the shortest odd cycle in $G$.)

We need one more observation, namely, that many of the components actually
contain only one vertex. Using our earlier notation,
let $v$ be an arbitrary vertex with neighbors $w_1,w_2,w_3$,
the corresponding random variables are $X$ and $Y_1,Y_2,Y_3$.
Note that $Y_1<0$, $Y_2<0$ and $Y_3<0$ imply that $X > 0$.
Therefore the probability $p$ that $v$ is an isolated vertex in
$G[V_{+}]$ is
%
%
\begin{eqnarray}\label{eqprobiso}
p &\stackrel{\mathrm{def}} {=}& P ( X > 0; Y_1<0;
Y_2<0; Y_3<0 ) = P ( Y_1<0;
Y_2<0; Y_3<0 ) \nonumber
\\
&=&
P ( y_i \cdot Z < 0; i=1,2,3 ) = \frac{1}{2} -
\frac{1}{4\pi} \sum_{1\leq i,j\leq3} \arccos
(c_{i,j})
\\
&=& \frac{1}{8} + \frac{1}{4\pi} \sum
_{1\leq i,j\leq3} \arcsin(c_{i,j}).\nonumber
\end{eqnarray}
Note that $\arcsin$ it is a monotone increasing
odd function on $[-1,1]$, which is convex on $[0,1]$.
Furthermore, the average of $c_{i,j}$ is $(\lambda_{\min}^2-3)/6 \geq
(2^2-3)/6 > 0$.
It is easy to see that these imply that
the right-hand side of (\ref{eqprobiso}) decreases (not increases)
if we replace each $c_{i,j}$ with their average $(\lambda_{\min
}^2-3)/6$. Thus
%
%
\begin{equation}
\label{eqprobisobound} p \geq\frac{1}{8} + \frac{3}{4\pi} \arcsin
\biggl(
\frac{\lambda
_{\min}^2-3}{6} \biggr).
\end{equation}
Our independent set will contain all isolated vertices
and at least a $(g-1)/(2g)$ proportion of all the other vertices in $V_{+}$.
This yields the following lower bound for the independence ratio of $G$:
\[
p + \frac{g-1}{2g} \biggl( \frac{1}{2} - p \biggr) =
\frac{g-1}{4g} + \frac{g+1}{2g} p.
\]
Combining this with (\ref{eqprobisobound}) yields the desired bound.

We can choose an independent set with the same expected size
from the ``negative vertices''
\[
V_{-} \stackrel{\mathrm{def}} {=} \bigl\{ v\in V(G)\dvtx
X_v < 0 \bigr\}.
\]
This implies the second part of the theorem.

We mention that the proof also works in the infinite setting,
so there is an analogous theorem for infinite transitive graphs
(as in Theorem~\ref{thminf}).
\end{pf}

%
\begin{remark}
Any nontrivial lower bound for the density of components
of size $3,5, \ldots$ in $G[V_+]$
would immediately yield an improvement in the above theorem.
In \cite{csghv} such nontrivial bounds
were obtained for the $3$-regular tree $T_3$.
\end{remark}
%

\begin{proposition} \label{proplamin3reg}
Suppose that $G$ is a finite, connected, $3$-regular, vertex-transitive graph.
Then either $G$ is isomorphic to the complete graph $K_4$,
or the least eigenvalue $\lambda_{\min}$ of its adjacency matrix is
at most $-2$.
\end{proposition}

The proof below is due to P\'eter Csikv\'ari.

\begin{pf*}{Proof of Proposition~\ref{proplamin3reg}}
Let $G$ be a connected, $3$-regular,
vertex-transitive graph with $\lambda_{\min}(G) > -2$.
We need to show that $G$ must be the complete graph $K_4$.

Cauchy's interlacing theorem implies that
$\lambda_{\min}(G) \leq\lambda_{\min}(H)$
whenever $H$ is an induced subgraph of $G$.
Therefore $\lambda_{\min}(H) > -2$ must hold for any induced subgraph.
Let $T$ denote the tree shown in Figure~\ref{figgr}.
It is easy to see that the smallest eigenvalue of $T$ is $-2$.
We also have $\lambda_{\min}(C_{2k}) = -2$ for the cycle of length
$2k$ for any $k \geq2$.
Therefore $G$ can contain neither $T$, nor $C_{2k}$ as an induced subgraph.
%
%
\begin{figure}[b]%

\includegraphics{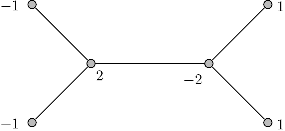}

%
%
%
%
\caption{The graph $T$ and the eigenvector corresponding to its least
eigenvalue $-2$.}
\label{figgr}
\end{figure}

We will distinguish three cases.
\begin{longlist}[\textit{Case}~3.]
\item[\textit{Case} 1.]
\textit{$G$ does not contain triangles.}

Let $u,v$ be two neighboring vertices, and let $u_1,u_2$ and $v_1,v_2$
denote the remaining two neighbors of $u$ and $v$, respectively.
Since $G$ contains no triangles, $u_1$, $u_2$, $v_1$, $v_2$ are pairwise
distinct vertices.
The induced subgraph on the set $\{u,u_1,u_2,v,v_1,v_2\}$ must be
isomorphic to $T$ (the graph shown in Figure~\ref{figgr}),
otherwise $G$ would contain a triangle or an induced $C_4$.
Since $G$ cannot contain $T$ as an induced subgraph, this is a contradiction.
\end{longlist}

\begin{longlist}[\textit{Case}~3.]
\item[\textit{Case} 2.]
\textit{$G$ contains triangles, but no two share a common edge.}

Since $G$ is vertex-transitive,
there must be at least one triangle through every vertex.
We claim that any two triangles must be disjoint.
If they had two common vertices,
then they would share an edge,
and if they had exactly one common vertex,
then that vertex would have degree at least $4$.

So we have disjoint triangles in $G$,
exactly one through every vertex.
We claim that there can be at most one edge between two triangles
(with one endpoint in one triangle and one in the other).
Indeed, otherwise we would either have an induced~$C_4$
or a vertex with degree at least $4$.

Let us consider the following graph $G^\ast$.
To each triangle in $G$ corresponds a vertex in $G^\ast$,
and we join two such vertices with an edge if
there is an edge between the corresponding triangles.
It is easy to see that $G^\ast$ will be $3$-regular as well.
Take a cycle in $G^\ast$ with minimum length $g \geq3$.
There is a corresponding cycle of length $2g$ in the original graph $G$.
It is easy to see that this must be an induced cycle, contradiction.
\end{longlist}

\begin{longlist}[\textit{Case}~3.]
\item[\textit{Case} 3.] \textit{$G$ contains two triangles sharing
an edge.}

Let $xy$ be an edge shared by triangles $xyu$ and $xyv$; see Figure~\ref{fig3}.

\begin{figure}

\includegraphics{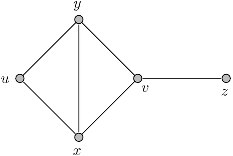}

\caption{Case 3.}\label{fig3}
\end{figure}

%
%
%

Then $x$ and $y$ already have degree $3$, while $u$ and $v$ still need
an edge.
We claim that $uv$ must be an edge.
Otherwise $v$ would have a neighbor $z$ different from $x,y,u$.
Since $z$ cannot be adjacent to $x$ and $y$,
there is only one triangle through $v$,
while there are two triangles through $x$,
contradticting the transitivity of $G$.
So $uv$ is an edge, and therefore each of $x,y,u,v$ has degree $3$.
Since $G$ is connected, $G$ cannot have any other vertices and thus
isomorphic to $K_4$.\quad\qed
\end{longlist}\noqed
\end{pf*}

For the sake of completeness we include a simple proof
for the following well-known result.
%

\begin{proposition} \label{propcorrdecay}
Let $G$ be a vertex-transitive graph, and
let $X_v$, $v \in V(G)$ be a factor of the i.i.d. process $Z_v$, $v
\in V(G)$
with $0< \var(X_v) < \infty$.
Then $\corr( X_v,X_{v'} ) \to0$ as the distance of $v$ and $v'$ goes
to infinity.
\end{proposition}

See \cite{bszv} for an explicit (and sharp) bound on the correlation
decay on the \mbox{$d$-}regular tree.

\begin{pf*}{Proof of Proposition \ref{propcorrdecay}}
For any vertex $v \in V(G)$, one can define the following L\'evy martingales:
\[
X_v^{(n)} = E\bigl( X_v \mid
Z_u, \dist(u,v) \leq n \bigr).
\]
According to martingale convergence theorems,
$X_v^{(n)}$ converges to $X_v$ almost surely and in $L^2$ as well.
Since $E(X_v^{(n)}) = E(X_v)$, the latter means that $\var
(X_v-X_v^{(n)}) \to0$ as $n \to\infty$.

Moreover, $X_v^{(n)}$, $v \in V(G)$ is a so-called \emph{block factor}
of $Z_u$, $u\in V(G)$,
that is, $X_v^{(n)}$ depends only on those $Z_u$'s
for which $u$ is in some finite neighborhood of $v$.

Now let $\varepsilon> 0 $ be arbitrary and let us pick $n$ such that
$ \var( X_v - X_v^{(n)} ) < \varepsilon$.
If the distance of $v$ and $v'$ is more than $2n$,
then $X_v^{(n)}$ and $X_{v'}^{(n)}$ are independent
(because they depend on disjoint sets of $Z_u$'s).
Therefore $\cov( X_v^{(n)},\break X_{v'}^{(n)} ) = 0$ and hence
\begin{eqnarray*}
\cov( X_v, X_{v'} ) &=& \cov\bigl( X_v-X_v^{(n)},
X_{v'}-X_{v'}^{(n)} \bigr) + \cov\bigl(
X_v^{(n)}, X_{v'}-X_{v'}^{(n)}
\bigr)
\\
&&{} + \cov\bigl( X_v-X_v^{(n)},
X_{v'}^{(n)} \bigr),
\end{eqnarray*}
which can be bounded by $\varepsilon+ 2 \sqrt{\varepsilon\var(X_v)}$,
and the statement of the proposition follows.
\end{pf*}

The following lemma is probably known,
but we did not find an explicit reference,
so we give a short proof.
%

\begin{lemma} \label{lemlama}
If $G$ is an infinite transitive graph, then the maximum $\lambda
_{\max}$ of the spectrum of $A_G$
is never in the point spectrum of $A_G$.
\end{lemma}

\begin{pf}
In the case $\lambda_{\max}=d$, the equation $A_G f= d f$ means that
the vector $f$ is harmonic.
However, the maximum principle implies that there are no $\ell^2$
harmonic functions.
Thus there is no eigenvector for $\lambda_{\max}$,
which is equivalent to saying that $\lambda_{\max}$ is not in the
point spectrum of $A_G$.

For the nonamenable case (i.e., $\lambda_{\max}<d$),
Theorem~II.7.8 in \cite{woessbook} implies that for any vertex $v$,
\[
\sum_{n=0}^\infty\lambda_{\max}^{-2n}
\bigl\langle1_v,A_G^{2n} 1_v
\bigr\rangle<\infty,
\]
where the left-hand side can be written in terms of the spectral
measure $\mu_G$ as
\[
\sum_{n=0}^\infty\lambda_{\max}^{-2n}
\int x^{2n} \,\mathrm{d}\mu_G(x) \ge\sum
_{n=0}^\infty\lambda_{\max}^{-2n}
\lambda_{\max}^{2n} \mu_G\bigl(\{
\lambda_{\max}\}\bigr).
\]
This forces $\mu_G(\{\lambda_{\max}\})=0$, which means that $\lambda
_{\max}$ is not in the point spectrum of $A_G$.
\end{pf}

%
\begin{lemma} \label{lemvolratio}
\[
\frac{\vol(S^{d-2})}{\vol(S^{d-1})} < \frac{\sqrt{d}}{\sqrt{2\pi
}}.
\]
\end{lemma}

\begin{pf}
Using the formula
\[
\vol\bigl( S^{n-1} \bigr) = \frac{2 \pi^{n/2}}{\Gamma(n/2)}, %
\]
we need to show that
\[
\frac{ \Gamma( (d-1)/2 ) }{ \Gamma( (d-2)/2 ) } < \sqrt{ \frac{d}{2} }.
\]
Since $\Gamma$ is log-convex, the increments of its logarithm over
intervals of length,
say, $1/2$ are increasing. Thus
\[
\frac{\Gamma((d-1)/2)}{\Gamma((d-2)/2)} \leq\frac{\Gamma({d}/2)}{\Gamma((d-1)/2)}
\]
and multiplying both sides by the left-hand side, we get
\[
\biggl(\frac{\Gamma((d-1)/2)}{\Gamma((d-2)/2)} \biggr)^2 \le
\frac{\Gamma(d/2)}{\Gamma((d-2)/2)}=
\frac{d-2}{2} < \frac{d}{2} %
\]
as required.
\end{pf}

%
\begin{lemma} \label{lemtangentline}
Let $\lambda\in[-3,-2]$ and
\[
f(t) \stackrel{\mathrm{def}} {=} \arcsin\biggl( \frac{ (1/2) -
\lambda+ ({\lambda^2}/{2}) + (1-\lambda) t }{
(1/2) - 2\lambda+ (\lambda^2/2) + t } \biggr).
\]
Then the tangent line to $f$ at $t_0 = (\lambda^2-3)/6$
is below $f$ on the entire interval $[-0.5,1]$; see Figure~\ref{figtangent} for the case $\lambda_{\min}=-2$.
\end{lemma}

\begin{pf}
We need to prove that
\[
f(t) - f'(t_0) t %
\]
takes its minimum value at $t_0$ on the interval $[-0.5,1]$.
This will follow from the fact that
$f'(t) < f'(t_0)$ for $-0.5 \leq t < t_0$ and
$f'(t) > f'(t_0)$ for $t_0 < t < 1$.

%
\begin{figure}

\includegraphics{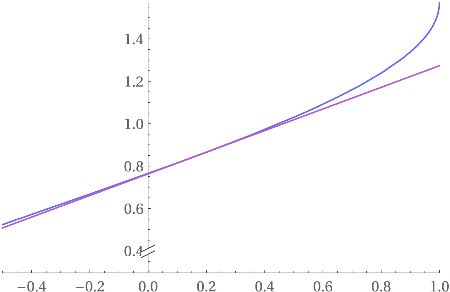}

\caption{The function $f(t)$ and its tangent at $t_0$ for $\lambda
_{\min}= -2$.}
\label{figtangent}
\end{figure}

In order to make calculations easier, we will use the following notation:
\begin{eqnarray*}
a &=& \frac{1}{2} - \lambda+ \frac{\lambda^2}{2} \geq4.5; \qquad b
= 1-
\lambda\geq3;
\\
c &=& a+b-1 \geq6.5;
\end{eqnarray*}
then
\[
f(t) = \arcsin\biggl( \frac{ a + b t }{ c + t } \biggr).
\]
It is easy to see that $ 0 < a+bt < c+t $ for $t \in[-0.5,1)$.
Therefore we have
\begin{eqnarray*}
f'(t) &=& \biggl( 1 - \biggl( \frac{ a + b t }{ c + t }
\biggr)^2 \biggr)^{-1/2} \frac{b(c+t) - (a+bt)}{(c+t)^2}
\\
&=&
\frac{bc - a}{(c+t) \sqrt{(c+t)^2 - (a+bt)^2}}.
\end{eqnarray*}
Since $bc-a > 0$ it follows that $f'$ is positive on $[-0.5, 1)$,
and thus $f$ is monotone increasing.
Next we study the intervals of monotonicity of $f'$.
First we note that
\[
(c+t)^2 - (a+bt)^2 = \bigl( c+a + (1+b)t \bigr)
\bigl(c-a+ (1-b)t \bigr).
\]
Using $c-a = b-1$ we get that
\[
(c+t)^2 - (a+bt)^2 = \bigl(b^2-1\bigr) (t+d)
(1-t),
\]
where
\[
d = \frac{c+a}{b+1} = \frac{1-3\lambda+\lambda^2}{2-\lambda} = 1-\lambda
- \frac{1}{2-\lambda} \geq
\frac{11}{4}.
\]
It follows that
\[
\frac{1}{(f'(t))^2} = \frac{b^2-1}{(bc-a)^2} (t+c)^2(t+d) (1-t).
\]
If we restrict ourselves to the interval $[-0.5,1)$
(where $f'$ is positive), then it suffices to examine the function
\[
g(t) = (t+c)^2(t + d) (1-t).
\]
Wherever $g$ is monotone increasing$,f'$ is monotone decreasing, and
vice versa.

So we have a fourth-degree polynomial $g$ with leading coefficient $-1$,
whose roots are $-c$ (with multiplicity $2$), $-d$ and $1$.
Consequently, the derivative $g'$ is a third-degree polynomial
with negative leading coefficient and
with roots $-c$, $u$, $v$, where $-c < u < -d < v < 1$.
We distinguish the following two cases.

\begin{longlist}[\textit{Case}~3.]
\item[\textit{Case} 1.]
$v \leq-0.5$.
Then $g$ is monotone decreasing on $[-0.5, \infty)$, and
therefore $f'$ is monotone increasing on $[-0.5,1)$,
and thus $f$ is convex on the whole interval,
which clearly implies the statement of the lemma.
\end{longlist}

\begin{longlist}[\textit{Case}~3.]
\item[\textit{Case} 2.]
$v > -0.5$.
Since the other two roots of $g'$ are less than $-d < -0.5$,
we know that $g$ is monotone increasing on $[-0.5,v]$ and
monotone decreasing on $[v,1)$. We claim that
%
%
\begin{equation}
\label{eqg} g \bigl( -\tfrac{1}{2} \bigr) > g \bigl(
\tfrac{1}{6} \bigr).
\end{equation}
This would yield that $v<1/6$.
Since $1/6 \leq t_0 = (\lambda^2-3)/6$, we have
$g(-1/2) > g(1/6) > g( t_0 )$.
This means that
$g(t) > g(t_0)$ for $-0.5 \leq t < t_0$ and
$g(t) < g(t_0)$ for $t_0 < t < 1$.
As for $f'$,
$f'(t) < f'(t_0)$ for $-0.5 \leq t < t_0$ and
$f'(t) > f'(t_0)$ for $t_0 < t < 1$,
and the statement of the lemma clearly follows.

It remains to show (\ref{eqg}).
Let $-1/2 = t_2 < t_1 = 1/6$.
Then $t_1 - t_2 = 2/3$; $t_2+c \geq6$ and $t_2+d \geq9/4$,
and consequently,
\begin{eqnarray*}
\frac{ g(t_1) }{ g(t_2) } &=& \frac{1-t_1}{1-t_2} \biggl( 1 + \frac
{t_1-t_2}{t_2+c}
\biggr)^2 \biggl( 1 + \frac{t_1-t_2}{t_2+d} \biggr)
\\
&\leq&
\frac{5/6}{3/2} \biggl( 1 + \frac{2/3}{6} \biggr)^2 \biggl(
1 + \frac{2/3}{9/4} \biggr) = \frac{5}{9} \biggl( \frac{10}{9}
\biggr)^2 \frac{35}{27} = \frac{17,500}{19,683} < 1.
\end{eqnarray*}
\end{longlist}\upqed
\end{pf}
\end{appendix}

\section*{Acknowledgments}
The authors are grateful to P\'eter Csikv\'ari
for the elegant proof of Proposition~\ref{proplamin3reg},
and to Gergely Ambrus, K\'aroly B\"or\"oczky, G\'abor Fejes T\'oth and
Endre Makai
for their remarks on Conjecture~\ref{conjgeom}.
We would also like to thank the anonymous referee
for the very careful reading of the manuscript
and the helpful comments and suggestions.




%

\printaddresses
\end{document}